\documentclass[12pt]{article}
\usepackage[cp1251]{inputenc}
\usepackage[russian, english]{babel}
\usepackage{amsfonts,amsmath}
\usepackage[labelsep=period]{caption}
\usepackage{tikz-cd}

\usepackage[left=1.5cm, right=1.5cm, top=2cm, bottom=2cm, headsep=0.7cm, footskip=1cm
]{geometry}

\newtheorem{theorem}{Theorem}
\newtheorem{lemma}{Lemma}
\newtheorem{corollary}{Corollary}
\newtheorem{remark}{Remark}
\newtheorem{example}{Example}

\begin{document}

\title{Integrable Birkhoff Billiards inside Cones}

\author {Andrey E. Mironov and Siyao Yin}
\date{}
\maketitle

\sloppy
\begin{abstract}
One of the most interesting problems in the theory of Birkhoff billiards is the problem of integrability. 
In all known examples of integrable billiards, the billiard tables are either conics, 
quadrics (closed ellipsoids as well as unclosed quadrics like paraboloids or cones), 
or specific configurations of conics or quadrics.
This leads to the natural question: are there other integrable billiards?
The Birkhoff conjecture states that if the billiard inside a convex, smooth, closed curve is integrable, then the curve is an ellipse or a circle.

In this paper we study the Birkhoff billiard inside a cone in $\mathbb{R}^n$.
We prove that the billiard always admits a first integral of degree two 
in the components of the velocity vector.
Using this fact, 
we prove that every trajectory inside a $C^3$ convex cone has a finite number of reflections.
Here, by $C^3$ convex cone, we mean a cone whose section with some hyperplane is a strictly convex closed $C^3$ submanifold of the hyperplane with nondegenerate second fundamental form.
The main result of this paper is the following. 
We prove that the Birkhoff billiard inside a convex $C^3$ cone is integrable.
This is the first example of an integrable billiard
where the billiard table is neither a quadric nor composed of pieces of quadrics.
\end{abstract}

\tableofcontents

\section{Introduction and Main Results}

The Birkhoff billiard is a dynamical system that studies the motion of a particle
in a domain $\Omega \subset \mathbb{R}^n$ with a piecewise smooth boundary $\partial \Omega$.
The particle moves freely with unit velocity vector inside the domain $\Omega$ and is elastically reflected upon colliding with the boundary: 
the angle of incidence equals the angle of reflection.
There are many remarkable results and unsolved problems in the theory of Birkhoff billiards (see e.g., \cite{Bir}--\cite{KT}). 
One of the most interesting problems is the Birkhoff conjecture (Birkhoff--Poritsky conjecture \cite{Por}).
Let $\Omega \subset \mathbb{R}^2$ be a convex domain and $\Sigma=\partial \Omega$ be a smooth closed curve (the {\it billiard table}). 
A curve $\Sigma_c \subset \Omega$ is called a {\it caustic} if every oriented line tangent to $\Sigma_c$ remains tangent to $\Sigma_c$ after the reflection. 
Lazutkin \cite{Laz} proved that 
if $\Sigma$ is sufficiently smooth with non-zero curvature, 
then there are infinitely many caustics near $\Sigma$. 
For example, if $\Sigma$ is an ellipse
$$
\frac{x^2}{a^2}+\frac{y^2}{b^2}=1,\quad a>b>0,
$$
then the confocal ellipse
$$
\frac{x^2}{a^2+\lambda}+\frac{y^2}{b^2+\lambda}=1,\quad a^2+\lambda>b^2+\lambda>0
$$
is a caustic for $\Sigma$. 
The domain bounded by the ellipse is foliated by caustics. 
The Birkhoff conjecture states that 
{\it if a neighborhood of the billiard table $\Sigma$ is foliated by caustics, then $\Sigma$ is an ellipse.}
There are many interesting partial results related to the Birkhoff conjecture (see e.g., \cite{B}--\cite{Kaloshin}), 
but in general it still remains open.
The Birkhoff billiard inside the ellipse is integrable; it admits a first integral $F$ (i.e., $F$ is constant along each trajectory), which is a polynomial in the components of the velocity vector $v = (v^1, v^2), |v| = 1$,
$$
F=a^2(v^2)^2+b^2(v^1)^2-(x(v^2)-y(v^1))^2.
$$
If the planar Birkhoff billiard admits a first polynomial integral in $v$, then the billiard table is an ellipse \cite{G} (see also \cite{Bol}, \cite{BM2}).

Another example of an integrable Birkhoff billiard 
(the Liouville--Arnold type integrability) is a billiard inside the ellipsoid
\begin{equation}\label{eq:ellipsoid}
\frac{(x^1)^2}{a_1^2} + \cdots + \frac{(x^n)^2}{a_n^2} = 1.
\end{equation}
There are $n-1$ polynomial first integrals $F_1(x,v),\dots,F_{n-1}(x,v)$ in involution
(see e.g., \cite{Ves})
$$
 \{F_i,F_j\}=\sum_{k=1}^n\left(\frac{\partial F_i}{\partial x^k}\frac{\partial F_j}{\partial v^k}-
 \frac{\partial F_i}{\partial v^k}\frac{\partial F_j}{\partial x^k}\right)=0.
$$
For example, at $n=3$ we have two integrals in involution
$$
F_1 = (v^1)^2 + \frac{(x^1 v^2- x^2 v^1 )^2}{a_1^2 - a_2^2}+ \frac{(x^1 v^3 - x^3 v^1)^2}{a_1^2 - a_3^2}, 
\quad 
F_2= (v^2)^2 + \frac{(x^1 v^2- x^2 v^1)^2}{a_2^2 - a_1^2}+ \frac{(x^2 v^3 - x^3 v^2)^2}{a_2^2 - a_3^2}.
$$
The integrability of Birkhoff billiards in various local senses has been investigated in \cite{Tr1}--\cite{G2}.

The concept of billiard caustics can be extended to higher dimensions.
The billiard inside the ellipsoid has confocal caustics defined by the equation of the form 
$$
\frac{(x^1)^2}{a_1^2+\lambda} + \cdots + \frac{(x^n)^2}{a_n^2+\lambda} = 1.
$$
Berger \cite{Ber1} proved that if $n \geq 3$ and the Birkhoff billiard inside a strictly convex $C^2$-smooth hypersurface $\Sigma \subset \mathbb{R}^n$ admits a strictly convex $C^2$-smooth caustic, then $\Sigma$ is an ellipsoid.

It turns out that, in the case of cones, the analogue of Berger's theorem is not true. 
There is always a family of caustics that are spheres. 
Our first result is the following.
\begin{figure}[htbp]
  \begin{center}
  \includegraphics[scale=0.20]{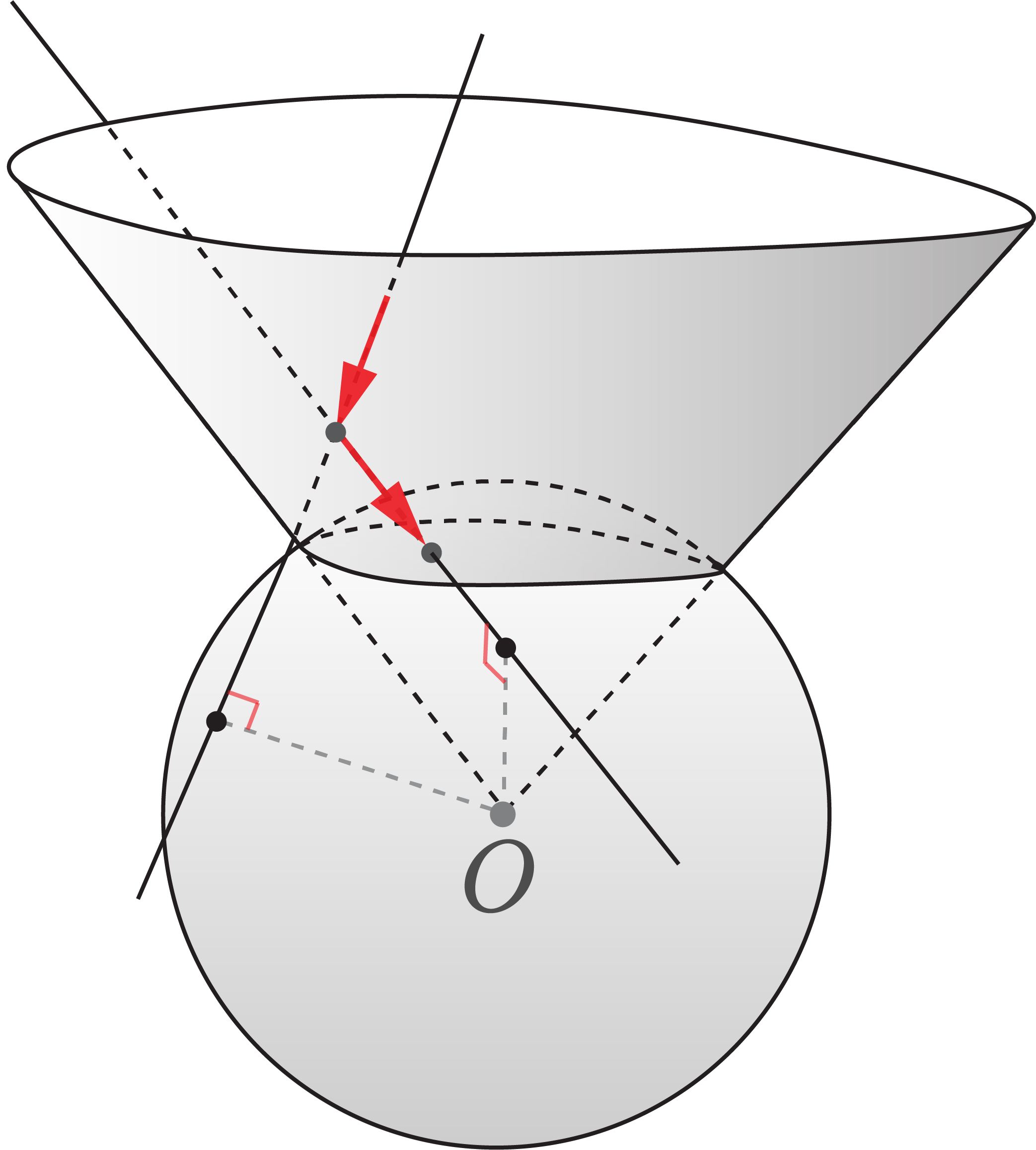}
  \end{center}
  \caption{The sphere as a caustic of the billiard inside a cone.}
  \label{fig:caustic}
\end{figure}
Let 
$\gamma \subset \{ x \in \mathbb{R}^n \mid x^n = 1 \}$ 
be a smooth
$(n-2)$-dimensional submanifold, 
and let 
$K = \{ tp \mid p \in \gamma, t > 0 \} \subset \mathbb{R}^n$ 
be an $(n-1)$-dimensional cone over $\gamma$. The following theorem holds.
\begin{theorem}\label{teo:cone-intro-n}
  1. 
  The Birkhoff billiard inside $K \subset \mathbb{R}^n$ admits the first integral
	$$
	I = \sum_{1 \leq i<j\leq n} m_{i,j}^2,
	$$
  where $m_{i,j} := x^i v^j - x^j v^i$ for $i < j$, $i, j = 1, \ldots, n$, and $v = (v^1, \ldots, v^n)$ is the velocity vector.
	
	2. 
  The spheres centered at the vertex $O \in \mathbb{R}^n$ of $K$ are caustics of the billiard inside $K$
  (see Fig. \ref{fig:caustic}).
\end{theorem}

\begin{figure}[h]
  \begin{center}
  \includegraphics[scale=0.4]{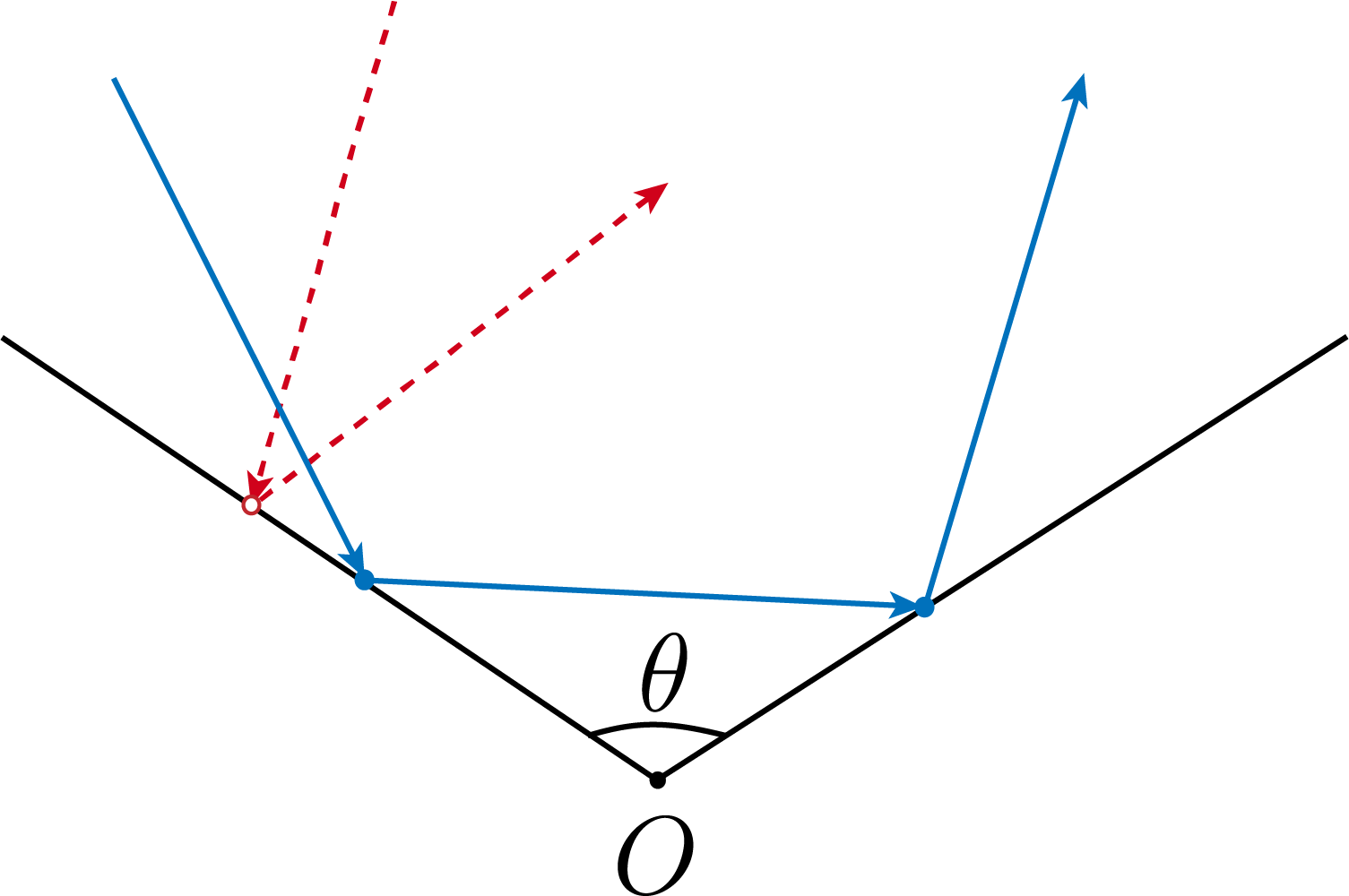}
  \end{center}
  \caption{For an angle $\theta$ in $\mathbb{R}^2$, a trajectory has $n$ or $n-1$ reflections, where $\frac{\pi}{n} \leq \theta < \frac{\pi}{n-1}$.
    }
  \label{fig:angle-theta}
  \end{figure}

Let us consider a billiard trajectory inside an angle (see Fig. \ref{fig:angle-theta}).
It is well known that the number of reflections is at most $\left\lceil \frac{\pi}{\theta} \right\rceil$,
where $\left\lceil x \right\rceil$ is the smallest integer greater than or equal to $x$ (see e.g., \cite{Tab}, \cite{Gal}).
Sinai \cite{Sin} studied Birkhoff billiards inside polyhedral angles in $\mathbb{R}^n$,
which are cones bounded by a set of hyperplanes.
It was proved that every trajectory has a finite number of reflections.
Moreover, there is a uniform estimate for the number of reflections for all trajectories in a fixed polyhedral angle (see also \cite{Sev}).

  In the next two theorems we assume that 
  $\gamma \subset \{ x \in \mathbb{R}^n \mid x^n = 1 \}$ 
  is a $C^3$-smooth closed submanifold whose second fundamental form is nondegenerate at every point. 
  From this it follows that $\gamma$ is strictly convex (see \cite{KN}).
  Let $K$ be a cone over $\gamma$. 
  Our next result is the following.

\begin{theorem}\label{teo:finite-intro}
  For any billiard trajectory inside $K$, the number of reflections is finite.
\end{theorem} 
\begin{remark}
  One can prove that there exist $C^2$-smooth convex cones that admit billiard trajectories with infinitely many reflections in finite time.
\end{remark}

Let us formulate our main result. 
The space of oriented lines in $\mathbb{R}^n$ admits a natural identification 
with the tangent bundle $TS^{n-1}$ of the unit sphere. 
Indeed, given an oriented line $l$ in $\mathbb{R}^n$, let $v \in S^{n-1}$ denote its unit direction vector. 
Let $Q \in l$ be the point that realizes the minimal distance to the origin, $\langle Q,v \rangle = 0$.
Hence $Q \in T_vS^{n-1}$. 
Therefore, each oriented line $l$ corresponds uniquely to a pair $(v,Q)$ in $TS^{n-1}$, establishing a bijection:
$$
l \longleftrightarrow (v,Q) \in TS^{n-1}.
$$

\begin{figure}[htbp]
	\centering
	\includegraphics[width=0.6\linewidth]{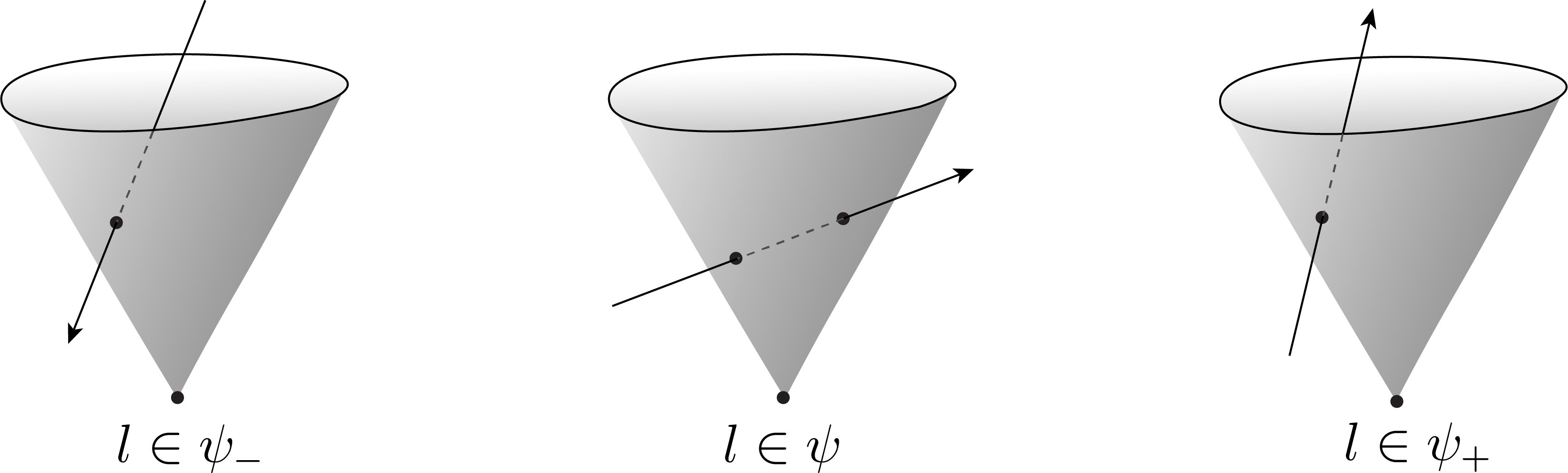}
  \caption{Oriented lines in $\psi_-$, $\psi$, and $\psi_+$, respectively. }
  \label{fig:cones-inte-s}
\end{figure}

The phase space $\Psi$ of the Birkhoff billiard inside $K$ is a subset of $TS^{n-1}$ consisting of oriented lines that intersect transversally with $K$. 
By Corollary \ref{lem:psi-open} (see Section 3.1 below) $\Psi$ is an open subset of $TS^{n-1}$, $\dim \Psi=2n-2$.
The phase space $\Psi$ naturally decomposes into three parts (see Fig. \ref{fig:cones-inte-s})
$$
\Psi = \psi_- \cup  \psi \cup \psi_+,
$$
where
$$
\begin{aligned}
	&\psi_- = \{\text{
		Oriented lines approaching from infinity and transversally intersecting the cone once
	}\}, \\
	&\psi \;\;= \{\text{
		Oriented lines intersecting  the cone twice
	}\}, \\
	&\psi_+ = \{\text{
		Oriented lines transversally intersecting the cone once before escaping to infinity
	}\}.
\end{aligned}
$$

We have the billiard map
$$
 \mu:\psi_-\cup\psi\rightarrow\psi\cup\psi_+,
$$
where $\mu$ maps an oriented line to the reflected oriented line.
Let us extend the billiard map on $\Psi$ by the identity map on $\psi_+$,
$$
 {\mu}(l)=l, \quad l\in\psi_+.
$$
Thus, we obtain the extended billiard map 
$$
{\mu}: \Psi \rightarrow \Psi. 
$$
Our main result is the following theorem.
\begin{theorem}\label{teo:integrability}
  There are continuous (smooth almost everywhere) first integrals 
  $I_1(x), \ldots, I_{2n-2}(x)$, $I_{2n-1}(x)=I$ on $\Psi$
invariant under $\mu$
  $$
  I_j(x)= I_j(\mu(x)), \quad  x \in  \Psi, 
  \quad j=1, \ldots, 2n-1.
  $$  
  Each point in the image of the map  
  $$
  \mathcal{I}= (I_1,\ldots, I_{2n-1}): \Psi \to \mathbb{R}^{2n-1}
  $$ 
  determines a unique billiard trajectory in $\Psi$.
\end{theorem}

There is a smooth submanifold $\Delta \subset \Psi$ (see Section 3) such that $I_j$ are smooth on $\Psi\setminus \Delta$. The first $2n-2$ integrals uniquely determine trajectories in $\Psi\setminus \Delta$, while the additional integral $I_{2n-1}$ is needed to distinguish trajectories in $\Delta$.

\begin{remark}
  We also construct smooth first integrals $I^s_1(x), \ldots, I^s_{2n-1}(x)$ on $\Psi$ 
  (see Lemma \ref{lem:many-integrals}, Section 3). The values of these integrals 
  uniquely determine billiard trajectories in $\Psi \setminus \Delta$ and vanish 
  identically on $\Delta$. Moreover, if we consider only the first $2n-2$ integrals, 
  whose number equals the dimension of $\Psi$, then at most two distinct trajectories 
  can be mapped to the same point in $\mathbb{R}^{2n-2}$ by 
  $(I^s_1(x), \ldots, I^s_{2n-2}(x))$.
  \end{remark}

  \begin{remark}
    The billiard inside the ellipsoid (\ref{eq:ellipsoid}), 
    as mentioned above, is integrable with $n-1$ first integrals $F_1, \ldots, F_{n-1}$ in involution. However, values of these integrals do not uniquely determine a billiard trajectory. In contrast, Theorem \ref{teo:integrability} provides integrability in a strong sense, where the values of the $2n-1$ integrals completely determine the trajectory.
  \end{remark}

The paper is organized as follows. 
We prove Theorems 1 and 2 in Sections 2, and Theorem 3 in Section 3.

The authors are grateful to Misha Bialy for valuable discussions and suggestions.

\section{Finite Number of Reflections inside $C^3$ Convex Cones}
In this section, we construct a first billiard integral, which is polynomial in the components of the velocity vector, and prove that every billiard trajectory inside a $C^3$ cone, under certain assumptions, has a finite number of reflections.

\subsection{Caustics and first integrals of degree two}
Here we prove Theorem \ref{teo:cone-intro-n}. First we prove the second part of the theorem, which follows immediately from the following lemma.

\begin{lemma}\label{lem:caus-dis}
	Let $K$ be a cone 
	in $\mathbb{R}^n$. 
	Suppose that the oriented line
	$l_1$ is reflected at 
	a point $p \in K$ to 
	the oriented line $l_2$. 
	Then the distance from $l_1$ 
	to the origin is equal to 
	the distance from $l_2$ 
	to the origin.	
\end{lemma}

\noindent\textbf{Proof.}
Let ${v}_1$ be the direction of $l_1$ 
and ${v}_2$ be the direction of $l_2$, 
with 
$\|{v}_1\| = \|{v}_2\| = 1$. 
Points on $l_i$ can be expressed as 
$p + t {v}_i$, $i=1,2$,
where $t \in \mathbb{R}$. 
To compute the distance from $l_i$ to the origin, 
we minimize 
$$
\|p + t {v}_i\|^2 = \|p\|^2 + 2t \langle p, {v}_i \rangle + t^2 \|{v}_i\|^2=
 \|p\|^2 + 2t \langle p, {v}_i \rangle + t^2.
$$
Solving 
$\frac{d}{dt} (\|p + t {v}_i\|^2) = 0$, we find $t = -\langle p,{v}_i \rangle$. 
Substituting back, we get
$$
\|p + t{v}_i\|^2_{\text{min}} = \|p\|^2 - \langle p,{v}_i \rangle^2.
$$
Thus
\begin{equation}\label{eqdistance}
  \text{dist}(l_i, O) = \left(\|p\|^2 - \langle p,{v}_i \rangle^2\right)^{1/2}.
\end{equation}
On the other hand, by the billiard reflection law,
$$
{v}_2 ={v}_1 - 2\langle{v}_1, {n} \rangle {n}.
$$
where $n$ is the normal vector to $K$ at point $p$.
From $\langle p, {n} \rangle = 0$ we obtain 
$$
\langle p,{v}_2 \rangle = \langle p,{v}_1 - 2\langle{v}_1, {n} \rangle {n} \rangle = \langle p,{v}_1 \rangle.
$$
Thus, from (\ref{eqdistance}) we have
$$
\text{dist}(l_1, O) = 
\text{dist}(l_2, O).
$$
Lemma \ref{lem:caus-dis} is proved.

\vspace*{1em}

The first part of Theorem \ref{teo:cone-intro-n} follows from the following lemma.

\begin{lemma}\label{lem:int-dis}
	For an oriented line $l \subset \mathbb{R}^n$
	with direction 
  ${v} = (v^1, \ldots, v^n)$, 
  $\|{v}\| = 1$, 
  the function 
	$$
  I = \sum_{1 \leq i<j \leq n} m_{i,j}^2 
  $$
	represents the square of the distance from
	$l$ to the origin, where
	$m_{i,j} := x^i v^j - x^j v^i$, $i < j$, $i,j = 1, \ldots, n$, 
	$(x^1, \ldots, x^n) \in l$.	
\end{lemma}
\textbf{Proof.}
Let $ {x} = (x^1, \ldots, x^n) $ be any point on $l$. Then by (\ref{eqdistance})
$$
\text{dist}(l, O)^2 = \|{x}\|^2 - \langle {x},{v} \rangle^2 
= \sum_{i=1}^n (x^i)^2 - \left( \sum_{i=1}^n v^i x^i \right)^2.
$$
Since $\sum_{i=1}^n (v^i)^2=1,$ this expands to
$$
\sum_{i=1}^n (x^i)^2 (1 - (v^i)^2) - 
\sum_{1\leq i < j \leq n} 2 v^i x^i v^j x^j 
=
\sum_{i=1}^n (x^i)^2\left(
 \sum_{j=1, j\neq i}^n
 (v^j)^2
  \right)-
\sum_{1\leq i < j \leq n} 2 v^i x^i v^j x^j=
$$
$$
 \sum_{1\leq i,j \leq n, i\neq j} (x^i)^2 (v^j)^2 
- \sum_{1\leq i < j \leq n} 2 v^i x^i v^j x^j
=
\sum_{1\leq i < j \leq n} (x^i v^j - x^j v^i)^2 
= \sum_{1\leq i < j \leq n} m_{i,j}^2 = I.
$$

Lemma \ref{lem:int-dis} is proved.

\subsection{Geometry of billiard trajectories inside cones}
In this section, we study the geometric properties of billiard trajectories inside a cone $K$ and prove Theorem \ref{teo:finite-intro}. 
Here $\gamma= K \cap \mathcal{P}$ is assumed to be a closed hypersurface of the hyperplane $\mathcal{P}=\{x\in\mathbb{R}^n\mid x^n=1\}$.
The proof of Theorem \ref{teo:finite-intro} is inspired by the results of Halpern \cite{Hal} and Gruber \cite{Gru2}. 
Halpern proved that in ${\mathbb R}^2$, $C^2$-smooth, strictly convex billiard tables can have trajectories with infinitely many reflections in finite time, whereas no such trajectories exist in the $C^3$ case. 
Gruber extended this result to ${\mathbb R}^n$, proving that $C^3$-smooth, strictly convex hypersurfaces do not admit trajectories with infinitely many reflections in finite time.

Consider a billiard trajectory
inside $K$ consisting of a sequence of 
oriented lines
$l_1, \ldots, l_n, \ldots$.
Let $l_k$ intersect $K$ at points 
$p_{k}$ and $p_{k+1}$, with direction vector ${v}_k$ of $l_k$ given by 
\begin{equation}\label{eq:vk}
  {v}_k = \frac{p_{k+1}-p_{k}}{\left\| p_{k+1}-p_{k}\right\|}.
\end{equation}
Let $n_k$ be the inward unit normal vector of $K$ at $p_k$.
The inward direction ensures 
$\langle {v}_k , n_k \rangle > 0$.

\begin{figure}[htbp]
  \begin{center}
  \includegraphics[scale=0.30]{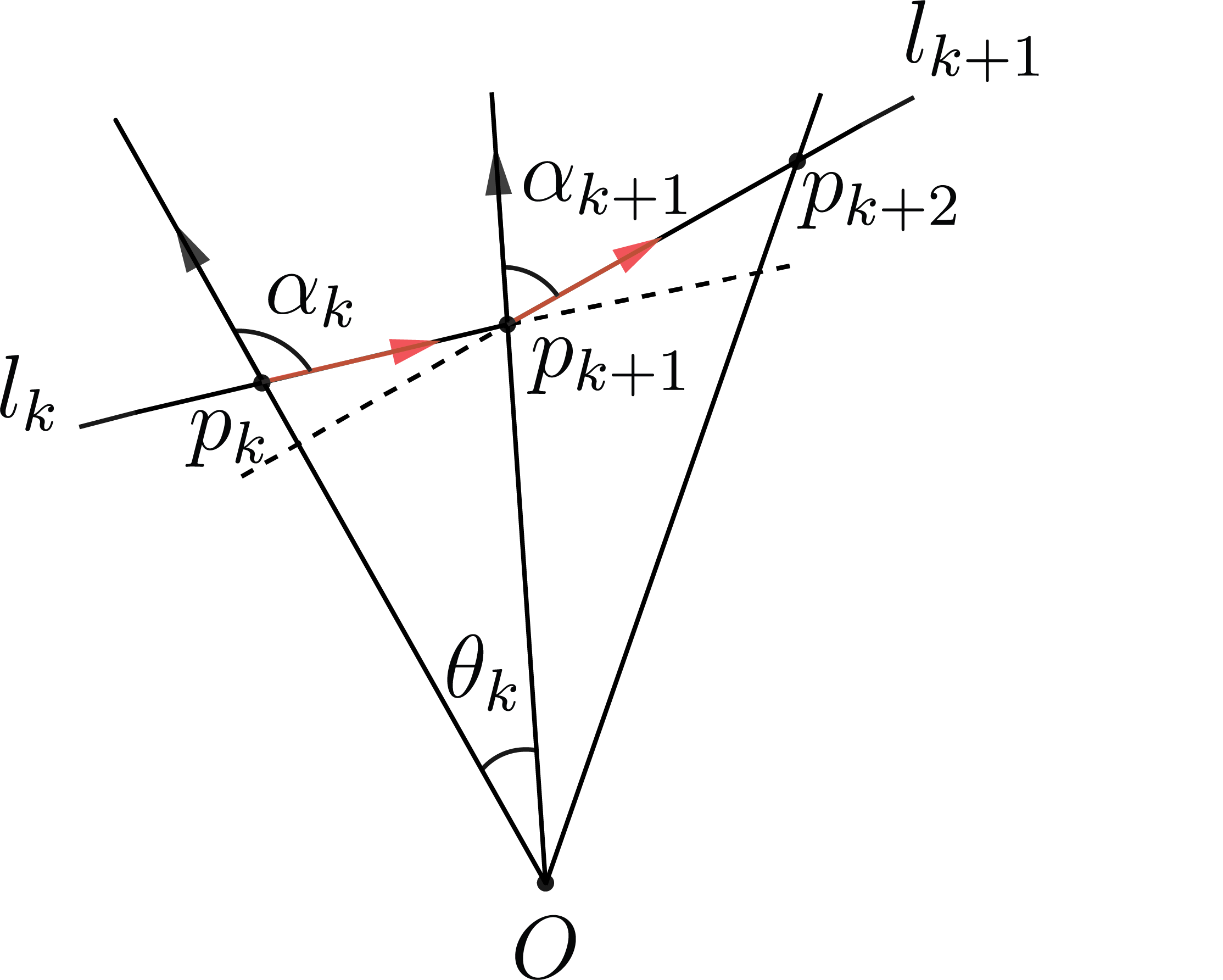}
  \end{center}
  \caption{The angles $\alpha_k, \theta_k$.}		
  \label{def-angles}
\end{figure}

We define angles $\alpha_k$, $\theta_k$ associated with the billiard trajectory (see Fig. \ref{def-angles}):
\begin{itemize}
  \item[1)] Let $\alpha_k$ be the angle between ${v}_k$ and the vector $p_k$, i.e.,
  $$
  \alpha_k = \arccos 
  \left(
    \left\langle
      {v}_k, \frac{p_k}{\|p_k\|} 
    \right\rangle
  \right)\in (0,\pi).
  $$

  \item[2)] Let $\theta_k$ be the angle $\angle p_k O p_{k+1}$,
            $$\theta_k \in (0,\pi).$$
\end{itemize}

In the next lemma we study properties of the angles $\alpha_k$, $\theta_k$, and points $p_k$. 
These properties do not require $K$ to be convex or $C^3$-smooth.

\begin{figure}[htbp]
	\centering
	\begin{minipage}[t]{0.33\linewidth}
		\centering
		\includegraphics[width=0.8\textwidth]{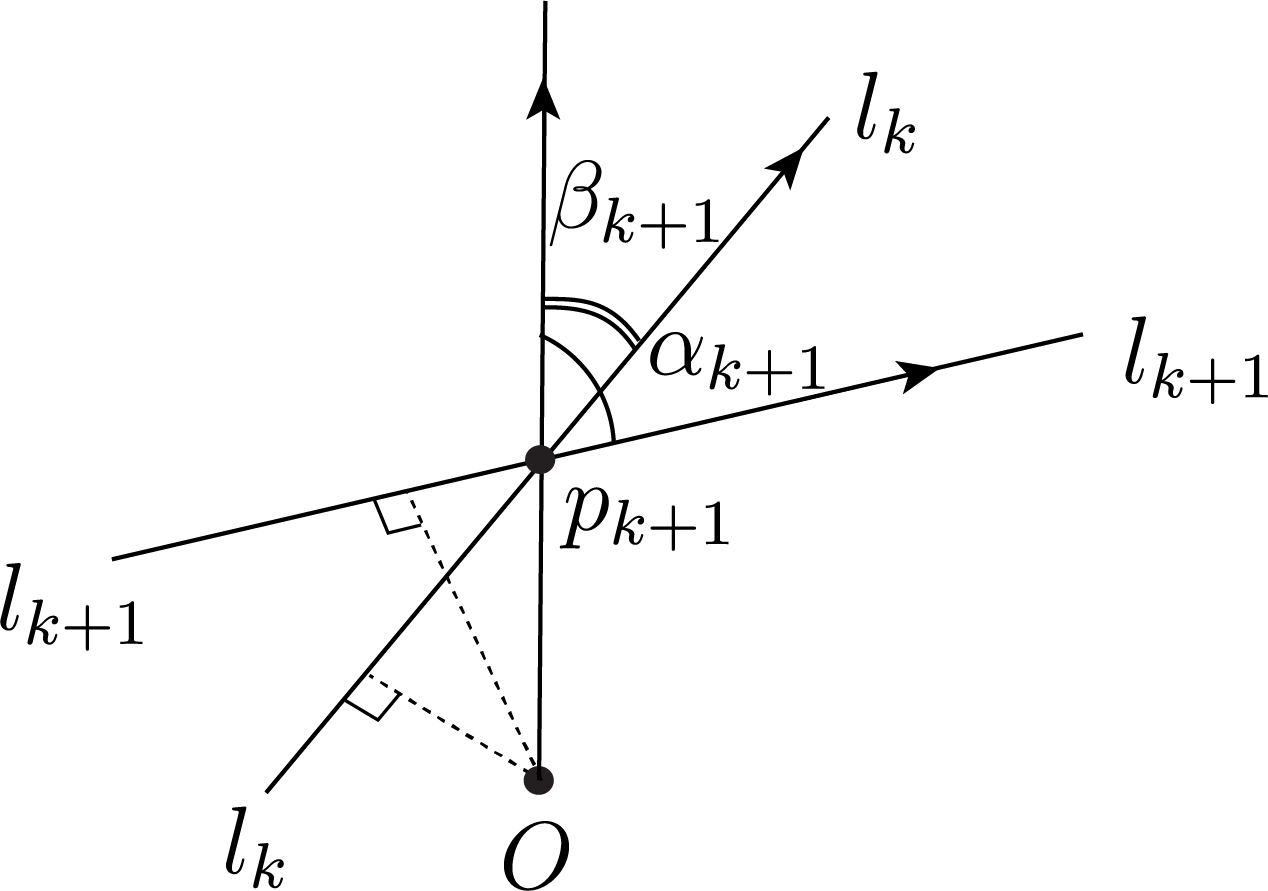}
    \caption{The angle $\beta_{k+1}$ between $l_{k}$ and $Op_{k+1}$. (Lines $l_k$, $l_{k+1}$ and $Op_{k+1}$ do not lie in the same plane.)}
    \label{fig:a-b1}
	\end{minipage}
  \hspace{1.6cm}
	\begin{minipage}[t]{0.32\linewidth}
		\centering
		\includegraphics[width=0.8\textwidth]{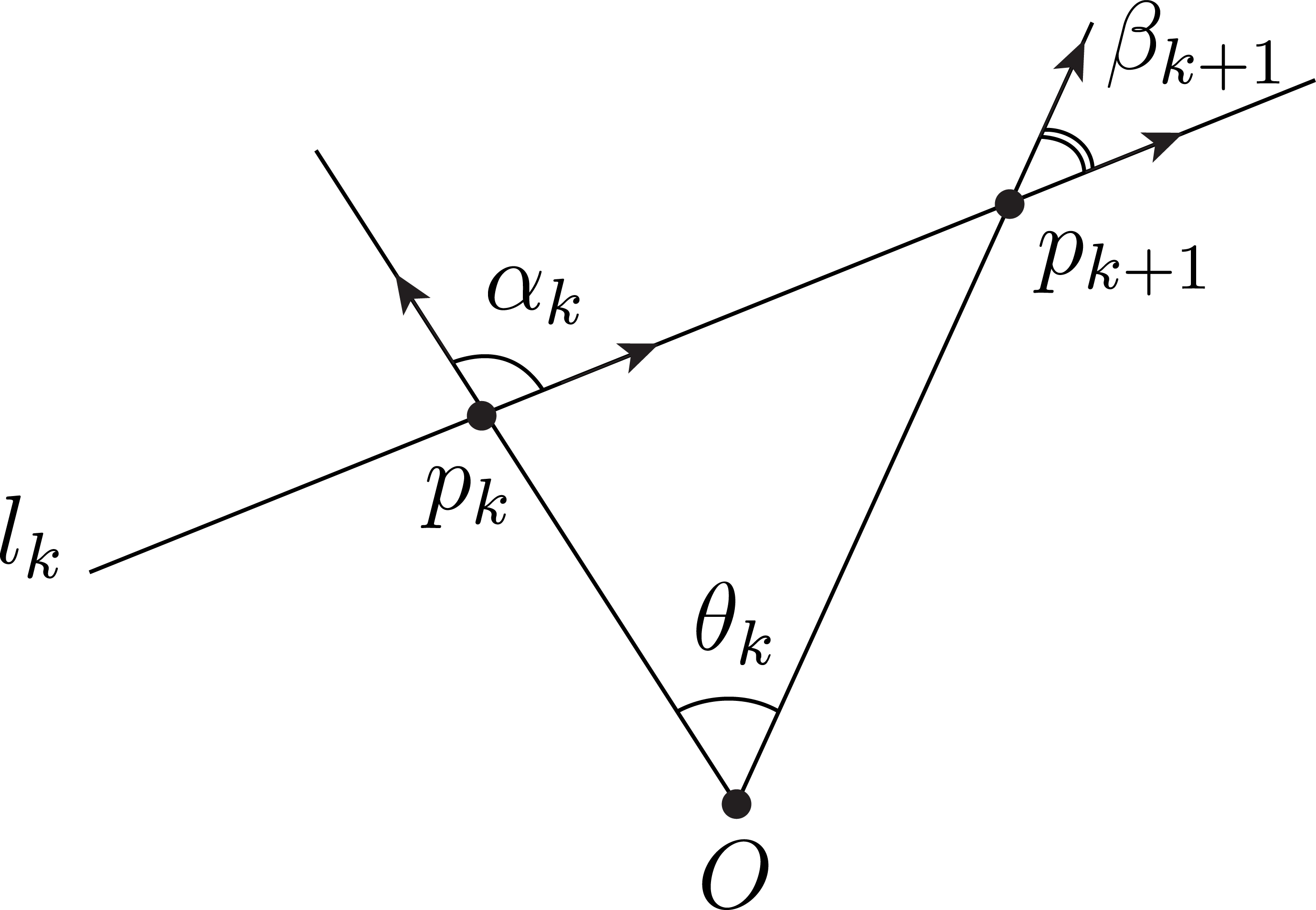}
    \caption{Correlation among $\alpha_k$, $\beta_{k+1}$ and $\theta_k$.}
    \label{fig:a-b2}
	\end{minipage}
  \end{figure}

\begin{lemma}\label{lem:theta-point}
    Suppose there exists a billiard trajectory inside $K$
    with infinitely many reflections. 
    Let $p_k = t_k q_k$, 
    where $t_{k} \in \mathbb{R}_{>0}$
    and $q_{k} \in \gamma$.
    Then the following statements hold:
    \begin{itemize}
    \item[1)] The angles $\alpha_k$ satisfy
    \begin{equation*}
      \alpha_{k+1} = \alpha_{k} - \theta_{k}.
    \end{equation*}
    Hence $\alpha_k$ decreases monotonically to $\alpha \geq 0$ as $k \to \infty$.
    
    \item[2)] The series $\sum_{k=1}^\infty \theta_k$ and 
    $\sum_{k=1}^\infty \| q_{k+1} - q_k \|$ converge.
    
    \item[3)] The points $q_k \to q \in \gamma$.
    
    \item[4)] If $\alpha > 0$, then $t_k \to t_0 > 0$ and $p_k \to p \in K$.
    If $\alpha = 0$, then $t_k \to +\infty$.
    \end{itemize}
  \end{lemma} 
\textbf{Proof.} 
1) 
Taking the inner product of $\frac{p_{k+1}}{\|p_{k+1}\|}$ 
with both sides of
$
{v}_{k+1} = {v}_k - 2\langle {v}_k, n_{k+1} \rangle n_{k+1} 
$  
and using $\langle n_{k+1}, p_{k+1} \rangle = 0$, we obtain
\begin{equation}\label{eq:vk-vkp1}
\left \langle {v}_{k+1}, \frac{p_{k+1}}{\|p_{k+1}\|}\right \rangle  
 = \left\langle {v}_{k}, \frac{p_{k+1}}{\|p_{k+1}\|} \right \rangle  
 - 2\left\langle {v}_k, n_{k+1} \right \rangle
 \left\langle n_{k+1}, \frac{p_{k+1}}{\|p_{k+1}\|} \right \rangle 
 = \left\langle {v}_k, \frac{p_{k+1}}{\|p_{k+1}\|}  \right\rangle .
\end{equation}
Let $\beta_{k+1}$ be the angle between line ${v}_k$ and $p_{k+1}$ (see Fig. \ref{fig:a-b1}). 
From (\ref{eq:vk-vkp1}), we have 
\begin{equation*}
  \cos \alpha_{k+1} = \cos \beta_{k+1},
\end{equation*} 
and since $\alpha_{k+1}, \beta_{k+1} \in (0, \pi)$, 
it follows that 
\begin{equation}\label{eq:a-b}
  \alpha_{k+1} = \beta_{k+1}.
\end{equation}

Next, consider $\triangle Op_k p_{k+1}$ 
(see Fig. \ref{fig:a-b2}). 
We have
\begin{equation}\label{eq:b-th}
\beta_{k+1} + \theta_k = \angle p_k p_{k+1} O + \theta_k = \pi - \angle p_{k+1} p_k O = \alpha_k.
\end{equation}
Combining 
(\ref{eq:a-b}) and (\ref{eq:b-th}), we obtain
$$
  \alpha_{k+1} = \alpha_{k} - \theta_{k}.
$$

Since $\theta_k > 0$, it follows that $\alpha_k$ is strictly decreasing. As $\alpha_k > 0$ for all $k$, we conclude that $\alpha_k \to \alpha \geq 0$ as $k \to \infty$.

  2) 
  From part 1) we know that
$$
 \sum_{k=1}^n\theta_k=\alpha_1-\alpha_{n+1},
$$
which implies that $\sum_{k=1}^{\infty} \theta_k $ converges to $\alpha_1 - \alpha$.

To show that $\sum_{k=1}^\infty \|q_{k+1} - q_{k}\|$ converges, we estimate $\|q_{k+1} - q_{k}\|$ in terms of $\theta_k$.
Let $\hat{l}_k$ be the line passing through $q_k$ and $q_{k+1}$.
The calculation of the area of $\triangle Oq_k q_{k+1}$ in two ways gives
\begin{equation}\label{eq:area-q}
\frac{1}{2}\|q_{k+1} - q_{k}\| \text{dist}(\hat{l}_{k}, O) =\frac{1}{2} \|q_{k+1}\| \|q_{k}\|\sin \theta_{k}.
\end{equation}
Since $q_k, q_{k+1}\in \gamma \subset \mathcal{P}$,
we have $\text{dist}(\hat{l}_{k}, O)\geq 1$.
Thus
$$
	\|q_{k+1} - q_{k}\| \leq
  \|q_{k+1}\|  \|q_{k}\|
  \sin \theta_k .
$$
Since $\gamma$ is bounded, there exists $a > 1$ such that $\|q \| \leq a$ for all $q \in \gamma$. Then
$$
\|q_{k+1} - q_{k}\| 
\leq a^2 \sin \theta_k
< a^2 \theta_k.
$$
Therefore $\sum_{k=1}^{\infty} \|q_{k+1} - q_{k}\|$ converges.

3)
From part 2), $\sum_{k=1}^\infty \|q_{k+1}- q_{k}\|$ converges, so $q_k$ is a Cauchy sequence in $\mathbb{R}^n$ and converges to a limit point $q$.
Since $\gamma$ is closed, the limit point $q \in \gamma$.

\begin{figure}[htbp]
	\centering
	\begin{minipage}[t]{0.22\linewidth}
		\centering
		\includegraphics[width=0.8\textwidth]{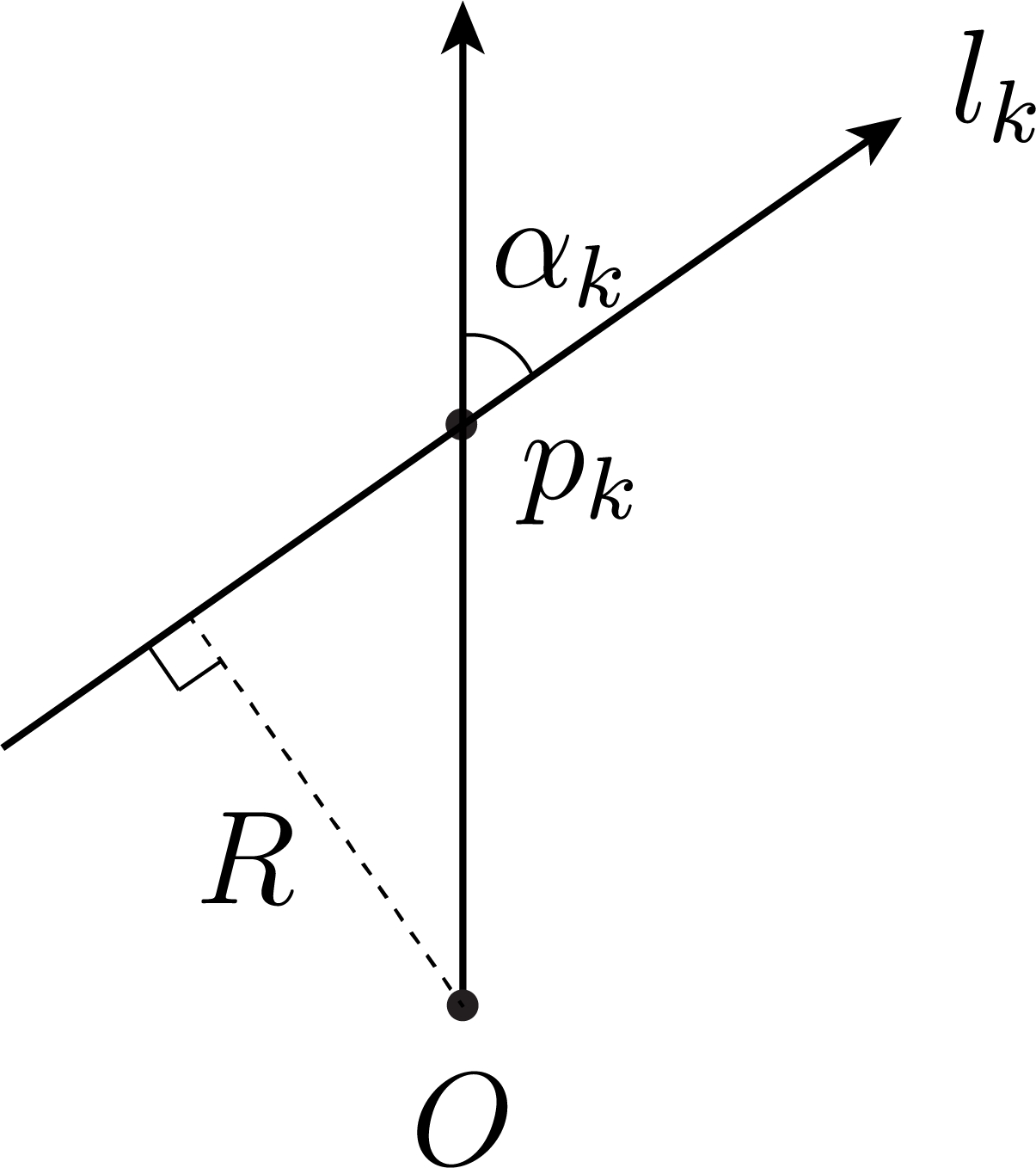}
		\caption{$\alpha_{k} < \frac{\pi}{2}$.}
    \label{fig:sin-alpha-1}
	\end{minipage}
  \hspace{3.5cm}
  \begin{minipage}[t]{0.20\linewidth}
		\centering
		\includegraphics[width=0.8\textwidth]{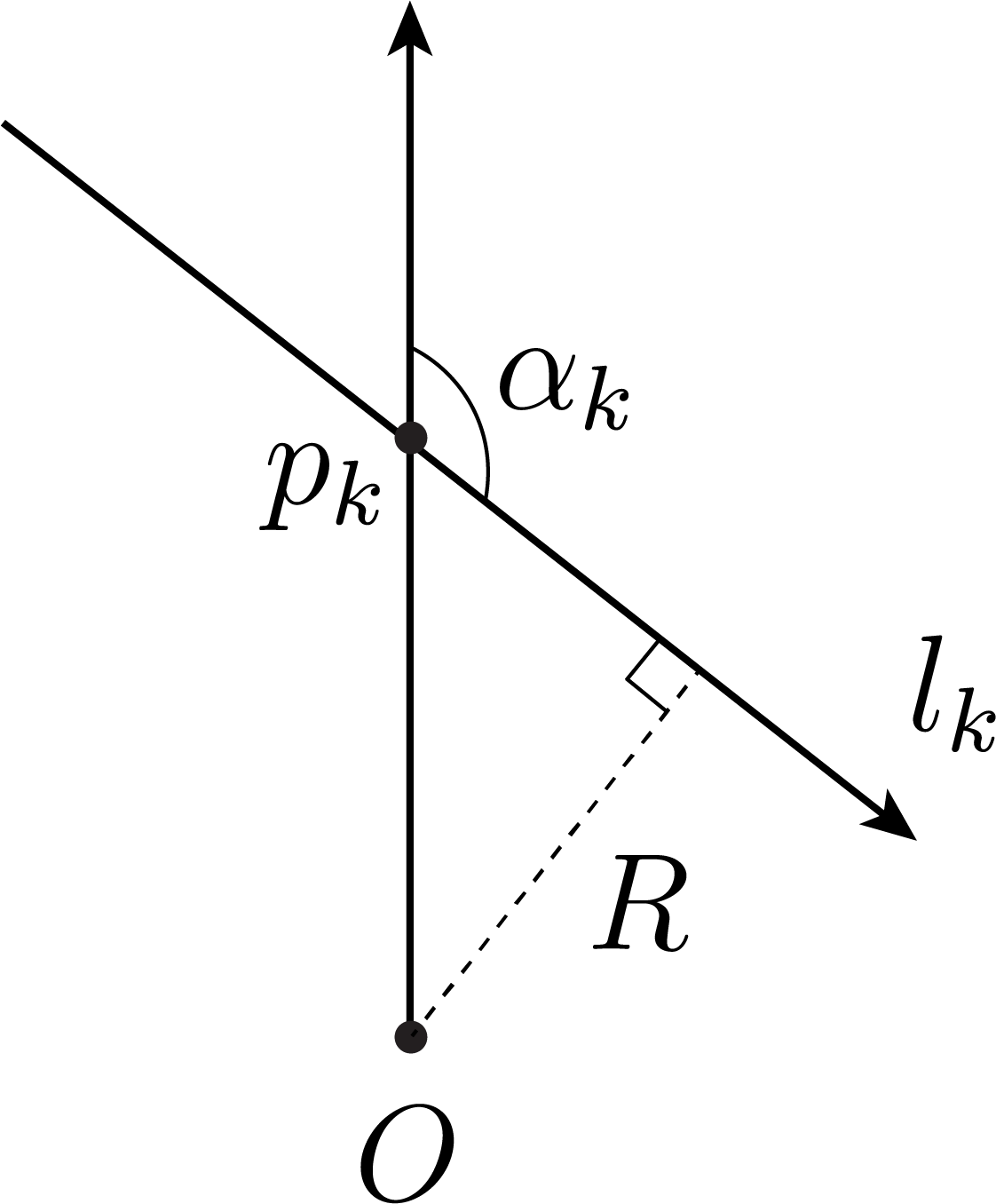}
		\caption{$\alpha_{k}>\frac{\pi}{2}$, \\ $\sin \alpha_k=\sin (\pi-\alpha_k)$.}
    \label{fig:sin-alpha-2}
	\end{minipage}
\end{figure}

4) 
Let $R$ be the distance from the origin to $l_k$. 
By Theorem \ref{teo:cone-intro-n}, 
$R$ is constant for all $k \geq 1$.
We have
\begin{equation}
\label{eq:sin-p}
\sin \alpha_k = \frac{R}{\|p_k\|}.
\end{equation}
This equality holds in all three possible cases: $\alpha_k < \frac{\pi}{2}$ (see Fig. \ref{fig:sin-alpha-1}); $\alpha_k > \frac{\pi}{2}$ (see Fig. \ref{fig:sin-alpha-2}, where we use $\sin \alpha_k = \sin(\pi - \alpha_k)$); and $\alpha_k = \frac{\pi}{2}$ (in this case, $R = \|p_k\|$).

From (\ref{eq:sin-p}), it follows
\begin{equation}\label{eq:sin-p1}
  \|p_k\| = \frac{R}{\sin \alpha_k},
\end{equation}
and then
\begin{equation}\label{eq:sin-t-q}
	t_k = \frac{\|p_k\|}{\|q_k\|}=\frac{R}{\sin \alpha_k \left\|q_k\right\|}.
\end{equation}
From (\ref{eq:sin-t-q}) and the convergence $\alpha_k \to \alpha$, $q_k \to q$, we obtain the following.
 If $\alpha = 0$, then $\lim_{k \to \infty} t_k = +\infty$. 
 If $\alpha > 0$, then $\lim_{k \to \infty} t_k = t_0$ for some $t_0 > 0$, 
and consequently $p_k \to p=t_0q \in K$.

Lemma  \ref{lem:theta-point} is proved.

\vspace*{1em}

 In the next lemma we study the properties of
 $\frac{\|p_k\|}{\|q_k\|}
 \langle  {v}_{k}, n_k \rangle$ 
 under the assumption that $K$ is $C^3$-smooth (i.e., $\gamma$ is $C^3$-smooth) and convex (i.e., $\gamma$ has everywhere non-degenerate second fundamental form).

  \begin{lemma}\label{lem:est}
    Let 
    $$
    A_k:=  \frac{\|p_k\|\langle  {v}_{k}, n_k \rangle}{\|q_k\|}, 
    \quad k\geq 1.
    $$
    Then:
    \begin{itemize}
      \item[1)] 
      There exist constants $M_1, M_2 > 0$ and an integer $k_0$ such that for all $k > k_0$
      \begin{equation*}
       M_1 \|q_{k+1} - q_k\|
       < A_{k}
       < M_2 \|q_{k+1} - q_k\|.
      \end{equation*}
      \item[2)]The difference satisfies  
      \begin{equation*}
       A_{k+1}
      -  
      A_{k}     
      =  O\left(
          A_k^2
        \right).
      \end{equation*}
    \end{itemize}
  \end{lemma}
\textbf{Proof.}
1)
Let us estimate $A_k$.
Since $ \langle p_k, n_k \rangle = 0$,
$p_{k+1}= \frac{\|p_{k+1}\|}{\|q_{k+1}\|} q_{k+1}$, 
we have
\begin{equation}\label{eq:ak-1}
 A_k = 
 \frac{\|p_k\|}{\|q_k\|} \left \langle \frac{  p_{k+1}-p_k}{ \|p_{k+1}-p_k \|}, n_k \right \rangle
=\frac{\|p_k\|\left \langle p_{k+1}, n_k \right \rangle}{\|q_k\|\|p_{k+1}-p_k \|} 
=
\frac{\|p_{k+1}\|\|p_k\|}{\|p_{k+1}-p_k\| \|q_{k+1}\|\|q_k\|}
\left\langle q_{k+1}, n_k  \right\rangle.
\end{equation}
The calculation of the area of $\triangle Op_k p_{k+1}$ in two ways gives
\begin{equation}\label{eq:area-p}
\frac{1}{2}\|p_{k+1} - p_{k}\| \text{dist}(l_{k}, O) =\frac{1}{2} \|p_{k+1}\|\|p_{k}\|\sin \theta_{k}.
\end{equation}
By Theorem \ref{teo:cone-intro-n}, 
$\text{dist}(l_{k}, O)$ is constant for all $k\geq 1$. 
Thus, we may assume $\text{dist}(l_{k}, O)=1$ without loss of generality. 
From (\ref{eq:area-q}) and (\ref{eq:area-p}), we obtain
\begin{equation}\label{eq:distance-p-q}
  \frac{\|p_{k+1}\|\|p_k\|}{ \| p_{k+1} - p_k\|\|q_{k+1}\|\|q_k\|} 
  = \frac{1}{\text{dist}(\hat{l}_k,O) \| q_{k+1} - q_k\|}.
\end{equation}
Substituting (\ref{eq:distance-p-q}) to (\ref{eq:ak-1})
we obtain
\begin{equation}\label{eq:ratio2}
A_k
= \frac{1}{\text{dist}(\hat{l}_k,O)}
\frac{ \left\langle q_{k+1}, n_k  \right\rangle}{\|q_{k+1}-q_k\|}.
\end{equation} 
Let $a$ be an upper bound for $\|q\|, q\in \gamma$. Then 
\begin{equation}\label{eq:bound-dist}
  1 \leq \text{dist}(\hat{l}_k,O) \leq \|q_k\|\leq a. 
\end{equation}
Since $n_k$ is the inward normal, we have $\langle {v}_{k}, n_k \rangle > 0$. Therefore $A_k > 0$ by definition, and from (\ref{eq:ratio2}), we know $\left\langle q_{k+1}, n_k \right\rangle > 0$. 
Combining this with (\ref{eq:bound-dist}), we have
\begin{equation}\label{eq:estimate-middle}
 \frac{1}{a} \frac{\left\langle q_{k+1}, n_k  \right\rangle}{\|q_{k+1}-q_k\|} 
 \leq 
 A_k 
 \leq 
 \frac{ \left\langle q_{k+1}, n_k  \right\rangle}{\|q_{k+1}-q_k\|}.
\end{equation}

To prove part 1) of Lemma \ref{lem:est}, it is enough to prove that 
there exist constants $M'_1, M'_2 > 0$ and an integer $k_0$ such that for all $k > k_0$
\begin{equation}\label{eq:estimate-middle3}
       M'_1 \frac{\left\langle q_{k+1}, n_k  \right\rangle}{\|q_{k+1}-q_k\|} 
       <
       \|q_{k+1}-q_k\| 
       < M'_2
       \frac{\left\langle q_{k+1}, n_k  \right\rangle}{\|q_{k+1}-q_k\|}.
\end{equation}
Indeed, from (\ref{eq:estimate-middle}) and (\ref{eq:estimate-middle3}) we obtain 
$$
\begin{aligned}
  & A_k 
  \leq  
  \frac{ \left\langle q_{k+1}, n_k  \right\rangle}{\|q_{k+1}-q_k\|}
  <
  \frac{1}{M_1'}  \|q_{k+1}-q_k\|,\\
  & A_k
  \geq 
  \frac{1}{a} \frac{\left\langle q_{k+1}, n_k  \right\rangle}{\|q_{k+1}-q_k\|} 
  >
  \frac{1}{aM_2'}  \|q_{k+1}-q_k\|,
\end{aligned}
$$
which proves part 1) with $M_1=\frac{1}{aM_2'}$ and $M_2=\frac{1}{M_1'}$.

Let us prove (\ref{eq:estimate-middle3}), 
which is equivalent to
\begin{equation}\label{eq:estimate-middle2}
  M'_1 \left\langle q_{k+1}, n_k  \right\rangle
  <
  \|q_{k+1}-q_k\|^2
  < M'_2
  \left\langle q_{k+1}, n_k  \right\rangle.
 \end{equation}

\begin{figure}[htbp]
  \begin{center}
  \includegraphics[scale=0.60]{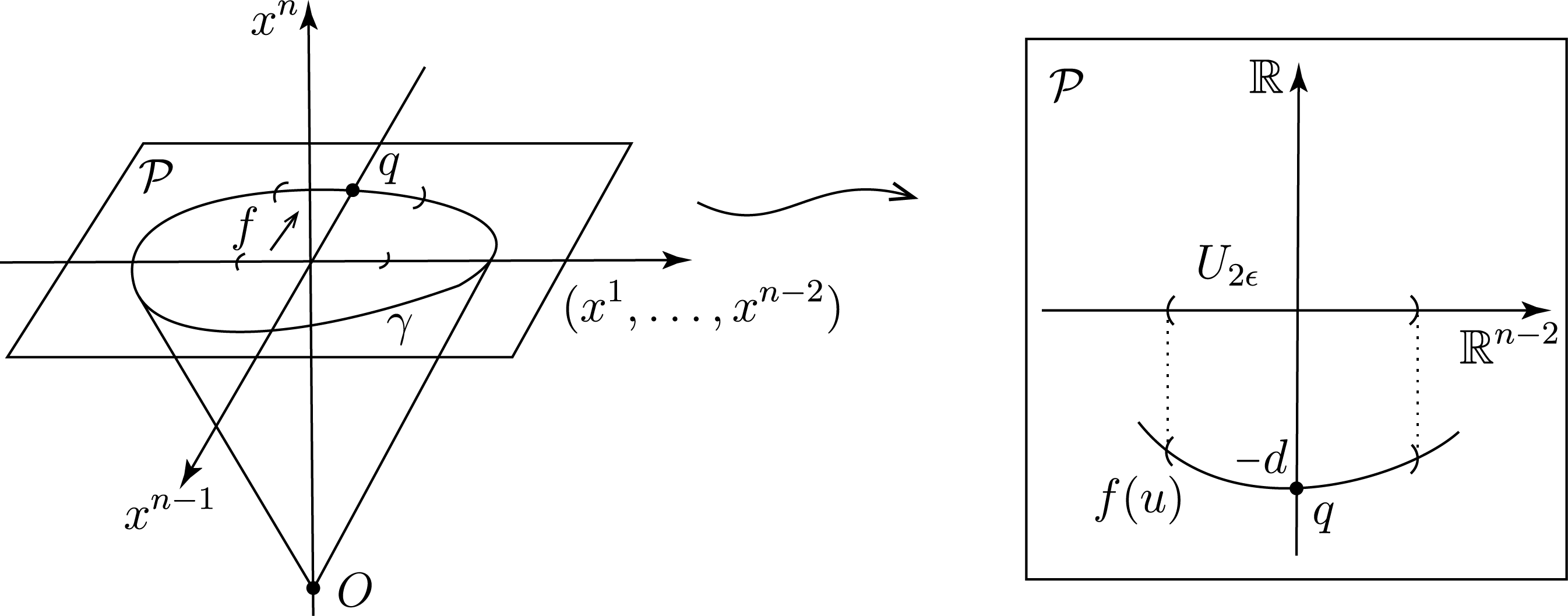}
  \end{center}
  \caption{The local coordinates near the limit point.}
  \label{fig:coordinates}
  \end{figure}

To obtain  (\ref{eq:estimate-middle2}), 
we introduce local coordinates near $q$. 
After a suitable rotation, we may assume $(0,\dots,0,1)$ lies inside $K$ and $q=(0,\dots,0,-d,1)$ for some $d>0$ 
(see Fig. \ref{fig:coordinates}).
Let $U_{\epsilon}$ be the open disc in $\mathbb{R}^{n-2}$ centered at 0 with radius $\epsilon$. Near $q$, since $\gamma$ is strictly convex,
we represent $\gamma$ as the graph of a $C^3$-smooth, strictly convex function $f: U_{2\epsilon} \to \mathbb{R}$ via
$$
\Phi(u):= (u, f(u), 1), \quad u= (u^1, \ldots, u^{n-2})\in U_{2\epsilon}.
$$ 
For large $k$, we can write $q_k = \Phi(u_k)$ with $u_k \in U_{\epsilon}$. 
Since $f$ is $C^3$ on $U_{2\epsilon}$, $f$ and the derivatives $f_i := \frac{\partial f}{\partial u^i}$, $f_{ij} := \frac{\partial^2 f}{\partial u^i\partial u^j}$ are bounded on the compact set $\overline{U}_{\epsilon}$.
Since $\gamma$ has non-degenerate second fundamental form everywhere, 
the Hessian matrix $\operatorname{Hess}f = (f_{ij})_{1\leq i,j\leq n}$
is positive definite on $\overline{U}_{\epsilon}$.

Let $V:= \{t\Phi(u)\mid u\in U_{2\epsilon}, t>0\} \subset K$.  Define an inward normal vector field on $V$
$$ 
n(u) := 
\left(
  -f_{1}(u),\ldots,-f_{{n-2}}(u), 1, -f(u)+ \sum_{i=1}^{n-2}u^if_{i}(u)
\right).
$$
This field is normal to $K$ since $\langle n(u), \Phi(u)\rangle = 0$, $\langle n(u), \frac{\partial \Phi}{\partial u^i}(u)\rangle = 0$ for $i=1,\dots,n-2$.
To verify $n(u)$ is inward, note that at $u=(0,\ldots,0)$, its last component equals $-f(0,\ldots,0)=d>0$. 
Since $(0,\ldots,0,1)\in \mathbb{R}^n$ lies inside $\gamma \subset \mathcal{P}$, this positive last component means $n(u)$ points to the inner side of the cone. As $n(u)$ is nowhere vanishing on $V$, it is an inward normal field on $V$.

To estimate $\left\langle q_{k+1}, n_{k} \right\rangle$
in (\ref{eq:estimate-middle2}), substituting
$q_{k+1}= \Phi( u_{k+1})$, 
$n_{k} = \frac{n(u_{k})}{\|n(u_{k})\|}$,
we have
\begin{equation}\label{eq:inner-product}
\left\langle q_{k+1}, n_{k} \right\rangle 
= 
\frac{ \left\langle \Phi( u_{k+1}), n( u_{k}) \right\rangle}{\|n(u_k)\|}.
\end{equation}
To evaluate $\left\langle \Phi( u_{k+1}), n( u_{k}) \right\rangle$, we have
\begin{equation}\label{eq:phi-n-1}
  \begin{split}
  \left\langle \Phi( u_{k+1}), n( u_{k}) \right\rangle
= & \left(-\sum_{i=1}^{n-2} u_{k+1}^i f_i(u_k )\right) 
+ f(u_{k+1}) - f(u_{k})
+ \sum_{i=1}^{n-2} u_{k}^i f_i(u_k ) \\
= &  f(u_{k+1}) -f(u_k ) 
- \sum_{i=1}^{n-2}(u_{k+1}^i - u_k^i) f_i(u_k ).
  \end{split}
\end{equation}
Since $f$ is $C^3$-smooth, from (\ref{eq:phi-n-1}) we have the Taylor expansion 
\begin{equation}\label{eq:convex-f}
  \left\langle \Phi( u_{k+1}), n( u_{k}) \right\rangle
= 
\frac{1}{2}
\left(
\sum_{i,j=1}^{n-2} f_{ij}(u_k )w^i w^j
\right)
\left\| u_{k+1}-u_k \right\|^2 + 
O(\left\|  u_{k+1} - u_k \right\|^3),
\end{equation}
where  $w=(w^1,\ldots,w^{n-2}):= \frac{u_{k+1}-u_k}{\left\|u_{k+1}-u_k\right\|}\in S^{n-3}$.
Since the Hessian matrix of $f$ is positive definite on $\overline{U}_{\epsilon}$ and $\overline{U}_{\epsilon}\times S^{n-3}$ is compact, 
there exists $c > 0$ such that
\begin{equation*}
\sum_{i,j=1}^{n-2} f_{ij}(u)w^i w^j \geq c,
\quad \forall (u,w)\in \overline{U}_{\epsilon}\times S^{n-3}.
\end{equation*}
Since $f$, $f_i$, $f_{ij}$ are bounded on $\overline{U}_{\epsilon}$,
there exists $C>0$ such that for all $u\in \overline{U}_{\epsilon} $ and $w\in S^{n-3}$
\begin{equation}\label{eq:bounded-n-fij}
1 \leq \|n(u)\|\leq C, 
\quad 
c\leq 
\sum_{i,j=1}^{n-2} f_{ij}(u)w^i w^j
\leq C.
\end{equation}
Hence, it follows from (\ref{eq:inner-product}), (\ref{eq:convex-f}), (\ref{eq:bounded-n-fij}) that 
\begin{equation}\label{eq:convex-f3}
  \left\langle q_{k+1}, n_k \right\rangle
= 
\frac{1}{2}
\frac{
\sum_{i,j=1}^{n-2} f_{ij}(u_k )w^i w^j}{\|n(u_k)\|}
\left\| u_{k+1}-u_k \right\|^2 + 
O(\left\|  u_{k+1} - u_k \right\|^3).
\end{equation}
In particular,
there exist constants $m_1, m_2>0$ and an integer $k_1$ such that for all $k > k_1$
\begin{equation}\label{eq:sin-u-bound}
m_1 \|u_{k+1}-u_k\|^2 
<
\left\langle q_{k+1}, n_{k} \right\rangle
< m_2 \|u_{k+1}-u_k\|^2.
\end{equation}

To estimate $\|q_{k+1} - q_{k}\|^2$
in (\ref{eq:estimate-middle2}), substituting 
$q_{k+1}= \Phi( u_{k+1})$, 
$q_{k}= \Phi( u_{k})$,
we have
\begin{equation}\label{eq:taylor-q-u1}
  \begin{split}
    \|q_{k+1} - q_{k}\|^2
   = &  \|\Phi(u_{k+1}) - \Phi(u_{k})\|^2 \\
    = &  \|u_{k+1}-u_k \|^2 
    +\left(f(u_{k+1})-f(u_k)\right)^2 \\
    = & \|u_{k+1}-u_k \|^2 
    + \left(\left(\sum_{i=1}^{n-2} f_{i}(u_k)w^i\right) 
    \|u_{k+1}-u_k \| 
    +O\left(\|u_{k+1}-u_k \|^2\right)\right)^2.
  \end{split}
\end{equation}
Since $\sum_{i=1}^{n-2} f_{i}(u)w^i$ is bounded for all  
$(u,w) \in \overline{U}_{\epsilon} \times S^{n-3}$, 
from (\ref{eq:taylor-q-u1}) we have
\begin{equation*}
    \|\Phi(u_{k+1}) - \Phi(u_{k})\|^2 
    = 
    \left(1 + \left(\sum_{i=1}^{n-2} f_{i}(u_k)w^i\right)^2  \right)
    \|u_{k+1}-u_k \|^2
    +O\left(\|u_{k+1}-u_k \|^3\right).
\end{equation*}
In particular,
there exist constants $m_3, m_4>0 $ and an integer $k_0 >k_1$ such that for all $k > k_0$
\begin{equation}\label{eq:q-u-bound}
m_3 \|u_{k+1}-u_k\|^2 
<
\|q_{k+1}-q_k\|^2 
< 
m_4 \|u_{k+1}-u_k\|^2.
\end{equation}

Combining (\ref{eq:sin-u-bound}) and (\ref{eq:q-u-bound}), we obtain for $k>k_0$, with $M'_1 = \frac{m_3}{m_2}$ and $M'_2 = \frac{m_4}{m_1}$,
\begin{equation*}
 M'_1 \left\langle q_{k+1}, n_k  \right\rangle
< \|q_{k+1}-q_k\|^2< 
M'_2\left\langle q_{k+1}, n_k  \right\rangle.
\end{equation*}
This proves (\ref{eq:estimate-middle2}) and hence part 1) of Lemma \ref{lem:est}.

2)  
Taking the inner product of ${n_{k+1}}$ with both sides of 
${v}_{k+1} = {v}_k - 2\langle {v}_k, n_{k+1} \rangle n_{k+1}$, 
we obtain
\begin{equation*}
  \langle  {v}_{k+1}, n_{k+1} \rangle
  = \langle  {v}_{k}, n_{k+1} \rangle
  -2 \langle  {v}_{k}, n_{k+1} \rangle
  \langle  n_{k+1}, n_{k+1} \rangle
  =- \langle  {v}_{k}, n_{k+1} \rangle.
\end{equation*}
Thus, using  $ \langle p_{k+1}, n_{k+1} \rangle = 0$,
$p_{k}= \frac{\|p_{k}\|}{\|q_{k}\|} q_{k}$, and (\ref{eq:distance-p-q})
\begin{equation*}
  A_{k+1}
  =\frac{\|p_{k+1}\|}{\|q_{k+1}\|}
  \left(-\langle  {v}_{k}, n_{k+1} \rangle\right)
  =\frac{\|p_{k+1}\|}{\|q_{k+1}\|}
  \left(-
  \left \langle \frac{  p_{k+1}-p_k}{ \|p_{k+1}-p_k \|}, n_{k+1} \right \rangle \right)
  =
  \frac{\|p_{k+1}\|}{\|p_{k+1}-p_k\| \|q_{k+1}\|}  
\left\langle p_{k}, n_{k+1}  \right\rangle
\end{equation*}
\begin{equation}\label{eq:ratio3}
=
\frac{\|p_{k+1}\|\|p_k\|}{\|p_{k+1}-p_k\| \|q_{k+1}\|\|q_k\|}  
\left\langle q_{k}, n_{k+1}  \right\rangle
= \frac{1}{\text{dist}(\hat{l}_k,O)}
\frac{\left\langle q_{k}, n_{k+1}  \right\rangle}
{\|q_{k+1}-q_{k}\|} 
\end{equation} 
To prove part 2) of Lemma \ref{lem:est}, by (\ref{eq:ratio2}) and (\ref{eq:ratio3}), it suffices to 
show that there exist a constant $M>0$ and an integer $k_2$ such that for all 
$k> k_2$ 
$$
\left|
A_{k+1} - A_k
\right|
=
\left|
\frac{1}{\text{dist}(\hat{l}_k,O)}
\frac{\left\langle q_{k}, n_{k+1}  \right\rangle}
{\|q_{k+1}-q_{k}\|} 
-
\frac{1}{\text{dist}(\hat{l}_k,O)}
\frac{\left\langle q_{k+1}, n_{k}  \right\rangle}
{\|q_{k+1}-q_{k}\|} 
\right| 
< M 
\left(
  \frac{1}{\text{dist}(\hat{l}_k,O)}
  \frac{\left\langle q_{k+1}, n_{k}  \right\rangle}
  {\|q_{k+1}-q_{k}\|} 
\right)^2,
$$
or equivalently,
\begin{equation}\label{eq:middle3}
\left|
\left\langle q_{k}, n_{k+1}  \right\rangle
-
\left\langle q_{k+1}, n_{k}  \right\rangle
\right| 
< M \frac{
\left\langle q_{k+1}, n_{k}  \right\rangle^2
}
{\text{dist}(\hat{l}_k,O)
\|q_{k+1}-q_{k}\|}.
\end{equation}

Now we prove (\ref{eq:middle3}).
To estimate the LHS of (\ref{eq:middle3}), we have 
\begin{equation}\label{eq:inner-product2}
\left\langle q_{k}, n_{k+1}  \right\rangle = 
\frac{\left\langle \Phi(u_{k}), n(u_{k+1}) \right\rangle}{\|n(u_{k+1})\|}.
\end{equation}
By the definition of $n(u)$ and $C^3$ smoothness of $f$ we have 
\begin{equation}\label{eq:convex-f5}
  \begin{split}
    \left\langle \Phi(u_{k}), n(u_{k+1}) \right\rangle
    = \ & \left(-\sum_{i=1}^{n-2} u_{k}^i f_i(u_{k+1})\right) 
    + f(u_{k}) - f(u_{k+1})
    + \sum_{i=1}^{n-2} u_{k+1}^i f_i(u_{k+1}) \\
    = \ & f(u_{k}) -f(u_{k+1}) 
    - \sum_{i=1}^{n-2}(u_{k}^i - u_{k+1}^i) f_i(u_{k+1}) \\
    = \ & \frac{1}{2}
    \left(
    \sum_{i,j=1}^{n-2} f_{ij}(u_{k+1})(-w)^i (-w)^j
    \right)
    \left\| u_{k+1}-u_k \right\|^2 + 
    O(\left\| u_{k+1} - u_k \right\|^3)\\
    = \ & \frac{1}{2}
    \left(
    \sum_{i,j=1}^{n-2} f_{ij}(u_{k+1})w^i w^j
    \right)
    \left\| u_{k+1}-u_k \right\|^2 + 
    O(\left\| u_{k+1} - u_k \right\|^3)\\
  = \ & \frac{1}{2}
    \left(
    \sum_{i,j=1}^{n-2}\left( f_{ij}(u_{k})+O(\left\| u_{k+1} - u_k \right\|)\right) w^i w^j
    \right)
    \left\| u_{k+1}-u_k \right\|^2 + 
    O(\left\| u_{k+1} - u_k \right\|^3)\\
    = \ & \frac{1}{2}
    \left(
      \sum_{i,j=1}^{n-2} f_{ij}(u_{k})w^i w^j
      \right)
      \left\| u_{k+1}-u_k \right\|^2 + 
      O(\left\| u_{k+1} - u_k \right\|^3).
  \end{split}
\end{equation}
We estimate the denominator of (\ref{eq:inner-product2}).
From (\ref{eq:bounded-n-fij}) we know $\|n(u)\|^{-1}$ is continuous on $\overline{U}_{\epsilon}$. Moreover, 
$\frac{\partial(\|n(u)\|^{-1})}{\partial u^k} = \frac{\text{a polynomial of } f, f_{i},f_{ij}, u^i}{\|n(u)\|^{3}}$ is also continuous on $\overline{U}_{\epsilon}$, thus 
$\|n(u)\|^{-1}$ is $C^1$ on $\overline{U}_{\epsilon}$.
By Taylor expansion, we have 
\begin{equation}\label{eq:norm}
  \|n(u_{k+1})\|^{-1}
  = \|n(u_{k})\|^{-1}
  +O(\|u_{k+1}-u_k\|).
\end{equation}
Hence, it follows from  (\ref{eq:bounded-n-fij}), (\ref{eq:inner-product2}), (\ref{eq:convex-f5}), and (\ref{eq:norm}) that 
\begin{equation}\label{eq:convex-f4}
  \left\langle q_{k}, n_{k+1} \right\rangle 
=
  \frac{1}{2} \frac{
    \sum_{i,j=1}^{n-2} f_{ij}(u_{k})w^i w^j
  }{\|n(u_{k})\|}
  \left\| u_{k+1}-u_k \right\|^2 + 
  O(\left\| u_{k+1} - u_k \right\|^3).
\end{equation}
Combining (\ref{eq:convex-f3}), (\ref{eq:convex-f4}) we obtain 
\begin{equation}\label{eq:lhs}
  \left\langle q_{k}, n_{k+1}  \right\rangle
-
\left\langle q_{k+1}, n_{k}  \right\rangle
= O(\|u_{k+1}-u_k \|^3).
\end{equation}
To estimate the RHS of (\ref{eq:middle3}), from (\ref{eq:bound-dist}), (\ref{eq:sin-u-bound}), and (\ref{eq:q-u-bound}), we have  for all $k > k_0$
\begin{equation}\label{eq:estimate5}
  \frac{
    \left\langle q_{k+1}, n_{k} \right\rangle^2
  }{
    \text{dist}(\hat{l}_k,O)
    \|q_{k+1}-q_{k}\|
  } 
  > \frac{m_1^2}{a} \frac{\|u_{k+1}-u_k\|^4}{\|q_{k+1}-q_{k}\|}
  > \frac{m_1^2}{a \sqrt{m_4}} \|u_{k+1}-u_k\|^3.
\end{equation}
Therefore, from (\ref{eq:lhs}), (\ref{eq:estimate5}), it follows that there exist a constant $M>0$ and an integer $k_2 > k_0$ such that for all 
$k> k_2$ 
$$
\left|
\left\langle q_{k}, n_{k+1}  \right\rangle
-
\left\langle q_{k+1}, n_{k}  \right\rangle
\right| 
< M 
\frac{m_1^2}{a \sqrt{m_4}} \|u_{k+1}-u_k\|^3
< M \frac{
\left\langle q_{k+1}, n_{k}  \right\rangle^2
}
{\text{dist}(\hat{l}_k,O)
\|q_{k+1}-q_{k}\|},
$$
i.e., (\ref{eq:middle3}) holds, 
which completes the proof of part 2) of Lemma \ref{lem:est}.

Lemma \ref{lem:est} is proved.

\subsection{Trajectories inside $C^3$ cones}
We are now ready to prove 
Theorem \ref{teo:finite-intro}. 
We need the following result by Weierstrass (see e.g., \cite{Kno}, p. 133).
\begin{lemma}[Weierstrass]\label{prop:w}
  Let $a_k$ be a series of non-zero reals such that 
  $$ \frac{a_{k+1}}{a_{k}} = 1 + b_{k},$$
  where $\sum_{k=1}^{\infty} |b_{k}|$ converges. 
  Then  $\sum_{k=1}^{\infty} a_k$ is divergent. 
\end{lemma}
Let us assume that there is a billiard trajectory inside $K$ consisting of an infinite series of reflections. Set
$$
a_{k} := A_k = \frac{\|p_k\|\langle  {v}_{k}, n_k \rangle}{\|q_k\|}.
$$
By Lemma \ref{lem:est}, part 1), we have 
$$
M_1 \|q_{k+1} - q_k\|
< a_{k}
< M_2 \|q_{k+1} - q_k\|, 
\quad
M_1, M_2>0.
$$
By Lemma \ref{lem:theta-point}, part 2), 
$\sum_{k=1}^{\infty} \left\|q_{k+1} - q_k\right\|$ converges,
hence $\sum_{k=1}^{\infty} a_k$ converges. 
On the other hand, 
by Lemma \ref{lem:est}, part 2), 
we have
\begin{equation*}
a_{k+1} - a_{k} = O(a_k^2),
\end{equation*}
hence 
\begin{equation*}
\frac{a_{k+1}}{a_k} - 1 = O(a_k).
\end{equation*}
Define 
$$
b_k:= \frac{a_{k+1}}{a_k} - 1,
$$
then $b_k = O(a_k)$. Since $a_k>0$ and $\sum_{k=1}^{\infty} a_k$ converges, $\sum_{k=1}^{\infty} |b_k|$ also converges, 
which implies that $\sum_{k=1}^{\infty} a_k$ diverges by Lemma \ref{prop:w}.
This contradiction shows that an infinite number of reflections cannot occur in a billiard trajectory inside $K$. 

This completes
the proof of Theorem \ref{teo:finite-intro}.
\begin{remark}
  We have used the condition that $K$ is $C^3$ (i.e., $f$ is $C^3$)
  in deriving (\ref{eq:convex-f3}), (\ref{eq:convex-f4}) in the proof of Lemma \ref{lem:est}, 
  which is crucial for the proof of the part 1 of Theorem \ref{teo:finite-intro}.
  Therefore, the method used in proving the part 1 of Theorem \ref{teo:finite-intro} does not apply to the $C^2$ case.
\end{remark}

\section{Integrability}
We establish the integrability of the billiard inside the $C^3$ convex cone $K$. 
Throughout this section, for any point $P \in K$, we denote by $n_P$ the inward unit normal vector to $K$ at $P$.

\subsection{The phase space}
In the beginning we prove that the phase space $\Psi$ as well as 
$\psi\cup\psi_+$ and $\psi\cup\psi_-$ are open in $TS^{n-1}$
(the space of all oriented lines in $\mathbb{R}^n$, see Introduction).

To prove this, we construct bijections from $\Lambda_1$ to $\psi \cup \psi_+$ 
and from $\Lambda_2$ to $\psi \cup \psi_-$,
where
$$
\Lambda_1 := \{(Q_1,v) \mid Q_1 \in K,  v \in S^{n-1}, \langle v, n_{Q_1} \rangle > 0\},
$$
$$
\Lambda_2 := \{(Q_2,v) \mid Q_2 \in K, v \in S^{n-1}, \langle v, n_{Q_2} \rangle < 0\}.
$$

\begin{figure}[htbp]
	\centering
	\begin{minipage}[t]{0.4\linewidth}
		\centering
		\includegraphics[width=0.8\textwidth]{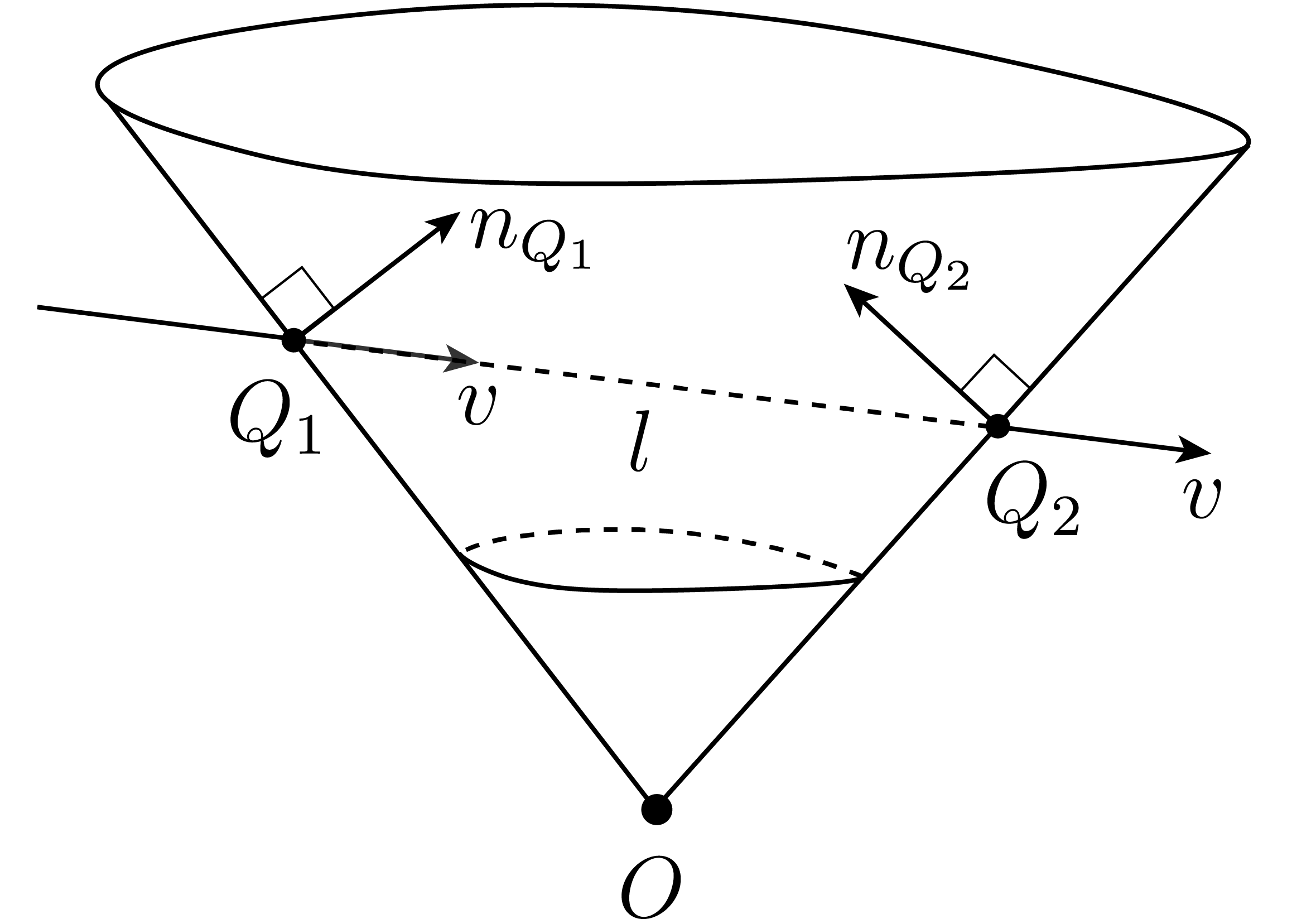}
    \caption{
    $\langle v, n_{Q_1} \rangle > 0$, $\langle v, n_{Q_2} \rangle < 0$,\\
    $l\in \psi$.
    }
    \label{fig:q1q2-psi}
	\end{minipage}
  \hspace{1cm}
  \begin{minipage}[t]{0.4\linewidth}
		\centering
		\includegraphics[width=0.8\textwidth]{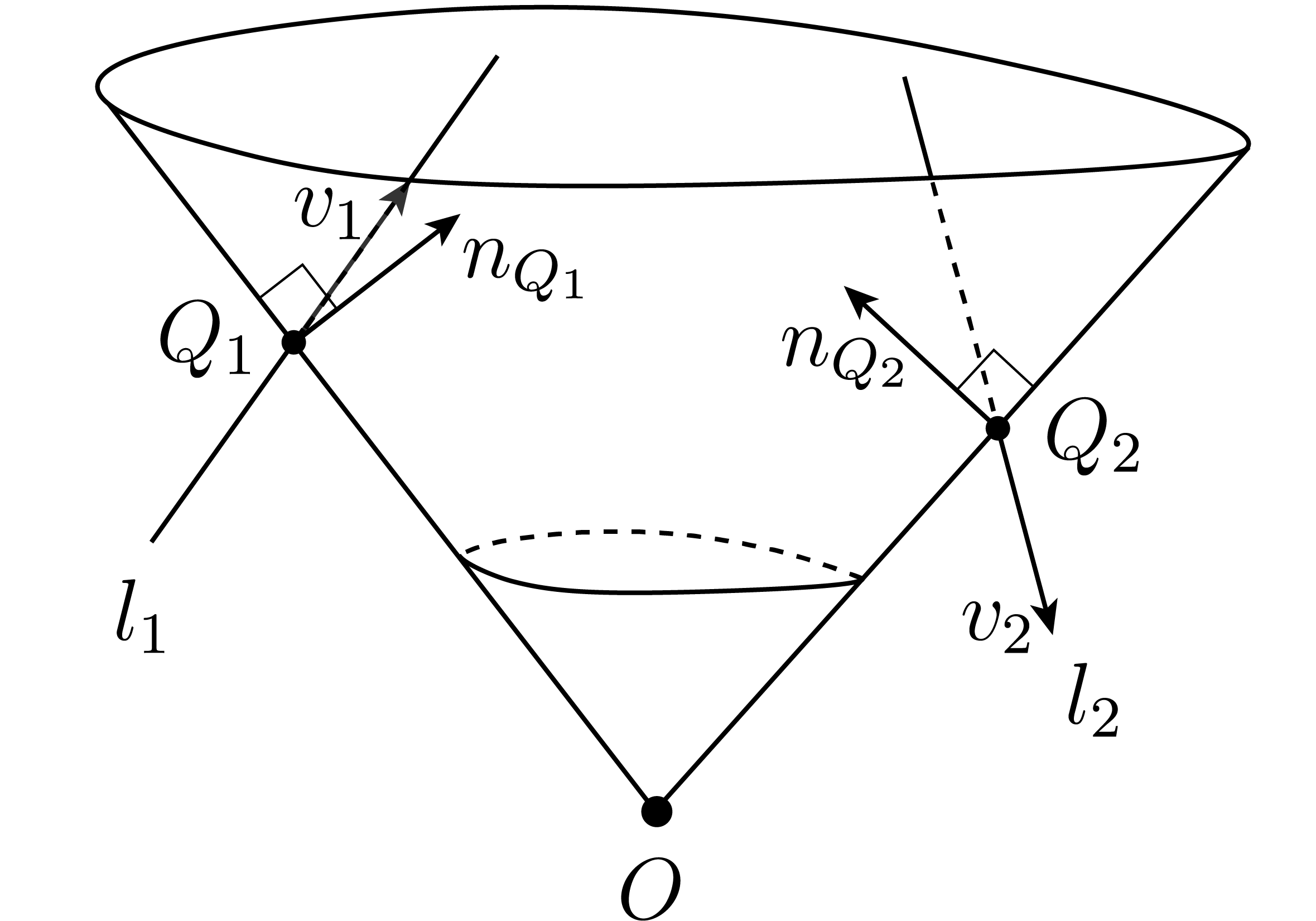}
    \caption{
    $\langle v_1, n_{Q_1} \rangle > 0$, $l_1 \in \psi_+$;\\
    $\langle v_2, n_{Q_2} \rangle < 0$, $l_2 \in \psi_-$.
    }
    \label{fig:q1q2-psi-pm}
	\end{minipage}
\end{figure}

For an oriented line $l$ in $\psi$, 
let $Q_1$ and $Q_2$ denote its first and second intersection points with $K$, respectively(see Fig. \ref{fig:q1q2-psi}).
For $l$ in $\psi_+$ (resp. $\psi_-$), 
which intersects $K$ at a single point,
let $Q_1$ (resp. $Q_2$) denote this intersection point 
(see Fig. \ref{fig:q1q2-psi-pm}).
For any oriented line $l$ in $\psi \cup \psi_+$, 
its direction vector $v$ has positive inner product with the inward normal vector at $Q_1$, 
hence corresponds to an element in $\Lambda_1$. 
Conversely, 
each pair $(Q_1,v)$ in $\Lambda_1$ uniquely determines an oriented line in $\psi \cup \psi_+$.
Similarly, 
there is a one-to-one correspondence between pairs $(Q_2,v)$ in $\Lambda_2$ 
and oriented lines in $\psi \cup \psi_-$, where the direction vector of each line has negative inner product with the normal vector at $Q_2$. Hence $\lambda_1$ and $\lambda_2$ are bijections.

Further we will use the following notations:
\begin{itemize}
  \item {\it $Q_1$, $Q_2$ are the intersection points of an oriented line with $K$ (the starting point and the end point, respectively) . }
  \item {\it $Q$ is the point on an oriented line that realizes the 
  distance between the oriented line and $O$.}
\end{itemize}

Since $\Lambda_1$ and $\Lambda_2$ are open subsets of $K \times S^{n-1}$,
they inherit smooth manifold structures from $K \times S^{n-1}$.
The bijections described above are explicitly given as follows, and they are diffeomorphisms.
\begin{lemma}\label{lem:valid-coords}
The mapping $\lambda_1:  \Lambda_1 \to \psi \cup \psi_+  \subset TS^{n-1}$,
$$
\lambda_1(Q_1,v)= (v,Q_1-\langle Q_1,v\rangle v),
$$
is a diffeomorphism.
  
Similarly, the mapping $\lambda_2:  \Lambda_2 \to \psi \cup \psi_-  \subset TS^{n-1}$,
$$
\lambda_2(Q_2,v) = (v,Q_2-\langle Q_2,v\rangle v),
$$
is a diffeomorphism.
\end{lemma}
\textbf{Proof.}
Let us check that $\lambda_1$ is well-defined, i.e., $(v,Q_1-\langle Q_1,v\rangle v) \in TS^{n-1}$.
We have $\langle Q_1-\langle Q_1,v\rangle v,v \rangle = \langle Q_1,v \rangle - \langle Q_1,v \rangle \langle v,v \rangle = 0$.

Since $\lambda_1$ is a bijection, to prove that 
$\lambda_1$ is a diffeomorphism, we need to show:
1) $\lambda_1$ is smooth, 
2) the differential $d\lambda_1$ is invertible at each point.
The smoothness is clear as $\lambda_1$ is composed of smooth operations.
To show $d\lambda_1$ is invertible, we compute its kernel explicitly.
Let $(Q_1(t),v(t))$ be a curve in $\Lambda_1$ with $(Q_1(0),v(0))=(Q_1,v)$.
Let $(\dot{Q_1},\dot{v})= \frac{{d}}{{d} t}(Q_1(t),v(t))|_{t=0}$ be the tangent vector of the curve at $t=0$.
Then
$$
\begin{aligned}
d\lambda_1(\dot{Q_1},\dot{v}) 
&=\frac{{d}}{{d} t} \big(v(t),Q_1(t)-\langle Q_1(t),v(t)\rangle v(t)\big)\big|_{t=0} \\
&=\big(\dot{v}, \dot{Q_1} - \langle \dot{Q_1},v \rangle v - \langle Q_1,\dot{v} \rangle v - \langle Q_1,v \rangle \dot{v}\big).
\end{aligned}
$$
For $d\lambda_1(\dot{Q_1},\dot{v})$ to be zero, we must have
$$
\dot{v} = 0 , \quad
\dot{Q_1} - \langle \dot{Q_1},v \rangle v = 0.
$$

The second equation implies $\dot{Q_1} = 0$. 
Indeed, the inner product of the second equation with $n_{Q_1}$ gives
$$ 
\langle \dot{Q_1},n_{Q_1} \rangle - \langle \dot{Q_1},v \rangle \langle v,n_{Q_1} \rangle = 0.
$$
Since $(Q_1,v) \in \Lambda_1$, we have $\langle v,n_{Q_1} \rangle > 0$. 
Moreover, as $\dot{Q_1}$ is tangent to $K$ at $Q_1$, we have $\langle \dot{Q_1},n_{Q_1} \rangle = 0$. 
These two conditions imply $\langle \dot{Q_1},v \rangle = 0$. 
Substituting this back into the second equation yields $\dot{Q_1} = 0$.

Thus, $d\lambda_1$ has trivial kernel. By the rank-nullity theorem and 
$\dim T_{(Q_1,v)}\Lambda_1 = \dim T_{\lambda(Q_1,v)}(TS^{n-1})$, $d\lambda_1$ is invertible at each point,
establishing that $\lambda_1$ is a local diffeomorphism to $TS^{n-1}$.
Since $\lambda_1$ is bijective between $\Lambda_1$ and $\psi \cup \psi_+$, 
it follows that $\lambda_1$ is a diffeomorphism from $\Lambda_1$ to $\psi \cup \psi_+$.

The proof for $\lambda_2$ is analogous to that for $\lambda_1$.

Lemma \ref{lem:valid-coords} is proved.
\vspace*{1em}

Since $\Lambda_{1}$, $\Lambda_{2}$ are open in $K\times S^{n-1}$,
the following corollary immediately follows from Lemma \ref{lem:valid-coords}.
\begin{corollary}\label{lem:psi-open}
  The sets $\psi \cup \psi_+$, $\psi\cup \psi_-$, and consequently their union $\Psi$ and intersection $\psi$, are open subsets of $TS^{n-1}$.
\end{corollary}

\vspace*{1em}

In the next lemma, 
we provide set-theoretic descriptions of $\psi_+$ and its boundary as subsets of $TS^{n-1}$. 
Let $\Gamma = K \cap S^{n-1}$.
$\Gamma$ bounds an open region $D \subset S^{n-1}$ inside $K$. 

The following lemma summarizes the convexity properties of $K$ that we will use (for a comprehensive introduction to convex geometry, see \cite{Sch}).

\begin{lemma}[Convexity Properties of $K$]\label{rem:convexity}
1)
For any $x, y \in K$, we have 
$\langle x,n_y \rangle \geq 0,$
where equality holds if and only if $x = \lambda y$ for some $\lambda > 0$. Consequently, 
\begin{equation}\label{eq:convexity1}
\langle x,n_y \rangle > 0 \iff \langle y,n_x \rangle > 0.
\end{equation}

2)
Let $\text{conv}(K)$ denote the convex hull of $K$. A point $p\in \mathbb{R}^n$ satisfies:
\begin{equation}\label{eq:convexity2}
 p \in \text{conv}{(K)}^{\circ} \iff  \langle p,n_v \rangle>0, \forall v\in \Gamma, 
\end{equation}
\begin{equation}\label{eq:convexity3}
p\in \mathbb{R}^n\setminus \overline{\text{conv}{(K)}} \iff  \exists v \in \Gamma, \langle p,n_v \rangle<0.
\end{equation}
Points satisfying the first condition (\ref{eq:convexity2}) are said to lie inside $K$;
those satisfying the second condition (\ref{eq:convexity3}) are said to lie outside  $K$.

\end{lemma}
\textbf{Proof.}
1)
Since $\gamma$ is a strictly convex hypersurface of $\mathcal{P} \cong
\mathbb{R}^{n-1}$, by the definition of $K$, we know that $K$ is a convex hypersurface of $\mathbb{R}^{n}$.
For any $y \in K$, $K$ is contained in the half-space $H_y^+$, where
$$
H_y^+ := \{z \in \mathbb{R}^n \mid \langle z,n_y \rangle \geq 0\}.
$$
Therefore,
$K \subset \cap_{y \in K} H_y^+$,
that is, for any $x, y \in K$, we have $\langle x,n_y \rangle \geq 0$.

If equality holds, i.e., there exist $x,y \in K$ such that $\langle x,n_y \rangle = 0$, then $x \in T_yK$. 
Let $x_1 = \frac{x}{x^n}, y_1 = \frac{y}{y^n} \in \gamma$,
where $x=(x^1,\ldots,x^n)$, $y=(y^1,\ldots,y^n)$.
We have $x_1 \in T_{y_1}\gamma= T_y K \cap \mathcal{P}$. 
By the strict convexity of $\gamma$, we have $x_1 = y_1$, 
hence $x = \frac{x^n}{y^n}y$, i.e., $x = \lambda y$ for some $\lambda > 0$.
Conversely, if there exists $\lambda > 0$ such that $x = \lambda y$, then clearly
$\langle x,n_y \rangle = \langle \lambda y,n_y \rangle = 0$.

The equivalence in \eqref{eq:convexity1} follows directly from the above conclusion.

2) 
First we prove \eqref{eq:convexity3}.
Since $\overline{\text{conv}{(K)}}$ is convex and closed as a subset of $\mathbb{R}^n$, we have (see e.g. Theorem 2.2.6 in \cite{Sch})
$$
\overline{\text{conv}{(K)}}= \cap_{y \in K} H_y^+.
$$
Then 
$$
 p\in \mathbb{R}^n\setminus \overline{\text{conv}{(K)}} 
 \iff  
 p\notin \cap_{y \in K} H_y^+
 \iff  
 \exists v \in \Gamma, \langle p,n_v \rangle<0.
$$

Next we prove \eqref{eq:convexity2}. 
Let us first prove that $p\in \text{conv}{(K)}^{\circ} \implies \langle p,n_v \rangle>0$ for all $v\in \Gamma$.
Let $p\in \text{conv}{(K)}^{\circ}$. By the definition of convex hull,
there exist a finite number of points $x_1,\ldots, x_k\in K$ and numbers $a_1,\ldots, a_k>0$, $\sum_{i=1}^{k}a_i=1$, such that 
$
p= \sum_{i=1}^{k}a_i x_i.
$
By part 1), for any $v \in \Gamma$, we have $\langle x_i,n_v\rangle \geq 0$ for all $i$. 
Therefore,
$$ 
\langle p,n_v\rangle = \sum_{i=1}^{k}a_i\langle x_i,n_v\rangle \geq 0 .
$$
The equality holds if and only if there exist $\lambda_i>0$ such that $x_i=\lambda_i v$ for all $i$.
This implies $p = (\sum_{i=1}^k a_i\lambda_i) v \in K$, contradicting $p\in \text{conv}{(K)}^{\circ}$.
Hence
$$ 
\langle p,n_v\rangle = \sum_{i=1}^{k}a_i\langle x_i,n_v\rangle > 0 .
$$
This proves that $p\in \text{conv}{(K)}^{\circ} \implies \langle p,n_v \rangle>0$ for all $v\in \Gamma$.

Conversely, suppose $\langle p,n_v \rangle>0$ for all $v\in \Gamma$. 
By \eqref{eq:convexity3}, we know that $p\in \overline{\text{conv}{(K)}}$.
We prove that $p$ must be in $\text{conv}{(K)}^{\circ}$.
Indeed, if not,
then $p\in K \cup\{O\}$.
There exists some $v\in \Gamma$ such that $\langle p,n_v \rangle=0$,
contradicting our assumption that $\langle p,n_v \rangle>0$ for all $v\in \Gamma$.

Lemma \ref{rem:convexity} is proved.

\vspace*{1em}

Recall that for a subset $A$ of a topological space $\mathcal{T}$, 
its boundary $\partial A$ is defined as 
$\partial A = \overline{A} \cap \overline{\mathcal{T} \setminus A}$.
In Lemma \ref{lem:boundary}, 
the notation $\partial \psi_+$ follows this definition.
The {\it boundary between $\psi_+$ and $\psi$} is defined as the subset of $\psi_+ \cup \psi$
consisting of their common limit points.

\begin{figure}[htbp]
  \begin{center}
  \includegraphics[scale=0.38]{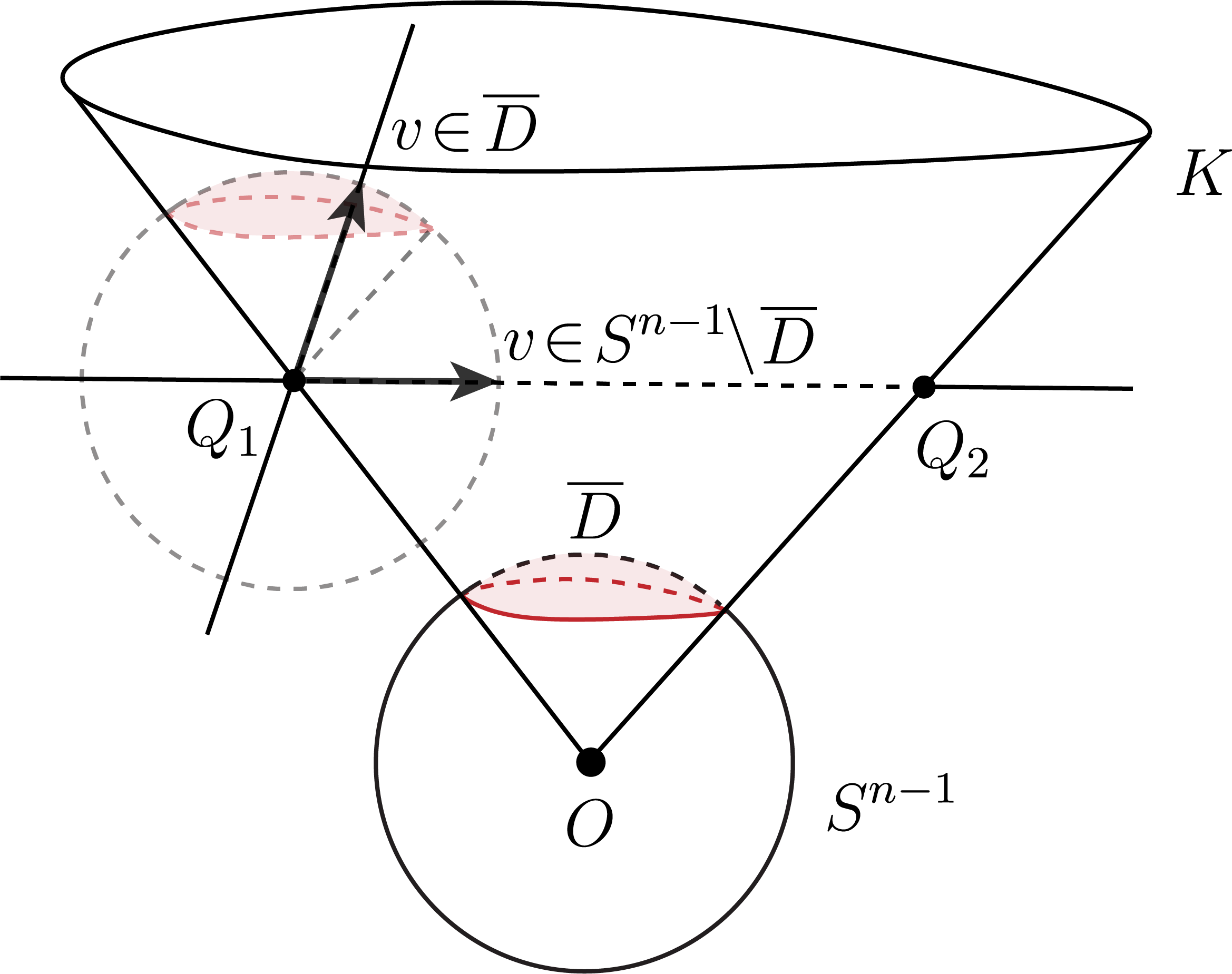}
  \end{center}
  \captionsetup{width=0.75\textwidth} 
  \caption{For $(Q_1,v) \in \Lambda_1$,  if $v \in \overline{D}$, $(Q_1,v)$ has no further intersection with $K$; 
  if $v \in S^{n-1}\setminus \overline{D}$, $(Q_1,v)$ intersects $K$ at a second point $Q_2$.}
  \label{fig:lambda}
\end{figure}

\begin{lemma}\label{lem:boundary}
  1) The boundary between $\psi_+$ and $\psi$ is
  $$
  \Delta_0 := \{(v,Q) \in TS^{n-1} \mid v \in \Gamma, \langle n_v,Q \rangle > 0\},
  $$
  which forms a codimension-1 submanifold of $\psi_+ \cup \psi$.
  
  2) $\psi_+$ admits the decomposition
  $$
  \psi_+ =\{(v,Q) \in TS^{n-1} \mid v \in D, Q\neq O\} \cup \Delta_0,
  $$
  where the first term represents its interior $\psi_+^{\circ}$.
  
  3) The boundary $\partial \psi_+$ is 
  $$
  \partial \psi_+ = \Delta_0 \cup \Delta_1 \cup \Delta_{\overline{D}},
  $$
  where   
  $$
  \Delta_1 := \{(v,Q) \in TS^{n-1} \mid v \in \Gamma, \langle n_v,Q \rangle \leq 0\},
  $$
  and
  $$
  \Delta_{\overline{D}} := \{(v,Q) \in TS^{n-1} \mid v \in \overline{D} , Q = O\}.
  $$
  \end{lemma}
\textbf{Proof.}
1) By Lemma \ref{lem:valid-coords}, we have a diffeomorphism $\lambda_1$ between $\Lambda_1$ and $\psi \cup \psi_+$. Geometrically, this diffeomorphism relates two equivalent representations of oriented lines: one using pairs $(Q_1,v)$ and another using pairs $(v,Q)$.
Let us characterize $\lambda_1^{-1}(\psi)$ and $\lambda_1^{-1}(\psi_+)$. 
Consider an oriented line $l$ intersecting $K$ at point $Q_1$ with $\langle v,n_{Q_1} \rangle>0$. The existence of a second intersection point $Q_2$ depends on the line direction: if $v \in \overline{D}$, $l$ has no further intersection with $K$; if $v \in S^{n-1}\setminus \overline{D}$, $l$ intersects $K$ at a second point $Q_2$ (see Fig. \ref{fig:lambda}). 

Thus we have
\begin{equation}\label{eq:lambda11}
\Lambda_1^1:=
\lambda_1^{-1}(\psi_+) =\{(Q_1,v) \mid Q_1 \in K, v \in \overline{D}, \langle v, n_{Q_1} \rangle > 0\},
\end{equation}
\begin{equation}\label{eq:lambda12}
\Lambda_1^2 :=
\lambda_1^{-1}(\psi) =\{(Q_1,v) \mid Q_1 \in K, v \in S^{n-1}\setminus \overline{D}, \langle v, n_{Q_1} \rangle > 0\}. 
\end{equation}
Let us define
\begin{equation}\label{eq:lambda10}
\Lambda_1^0 := \{(Q_1,v) \mid Q_1 \in K, v \in \Gamma, \langle v, n_{Q_1} \rangle > 0\}\subset \Lambda_1^1.
\end{equation}
From the set descriptions we see that $\Lambda_1^0$ is the boundary between $\Lambda_1^1\cong \psi_+$ and $\Lambda_1^2\cong \psi$, 
and is a codimension-1 submanifold of $\Lambda_1$.
Thus, to prove 1) it suffices to show that
\begin{equation}\label{eq:lambda-delta}
\lambda_1(\Lambda_1^0) = \Delta_0.
\end{equation}

First, we show that $\lambda_1(\Lambda_1^0) \subset \Delta_0$.
For any $(Q_1,v) \in \Lambda_1^0$, let $(\tilde{v},Q) = \lambda_1(Q_1,v)$.
By definition of $\lambda_1$, $\tilde{v} = v$ and $Q = Q_1-\langle Q_1,v\rangle v$.
Since $(Q_1,v) \in \Lambda_1^0$, we have $v \in \Gamma$ and $\langle v,n_{Q_1} \rangle > 0$.
Thus, from (\ref{eq:convexity1}) in Lemma \ref{rem:convexity}, part 1), we have
$$
\langle n_v,Q_1 \rangle > 0.
$$
Then,
$$
\langle n_v,Q \rangle = \langle n_v,Q_1-\langle Q_1,v\rangle v \rangle = \langle n_v,Q_1 \rangle > 0.
$$
Therefore $\lambda_1(Q_1,v) \in \Delta_0$.

To prove the reverse inclusion $\Delta_0 \subset \lambda_1(\Lambda_1^0)$, for any $(v,Q) \in \Delta_0$ it is necessary to prove that the oriented line $(v,Q)$ intersects $K$ at some point $Q_1$ where $\langle v, n_{Q_1}\rangle>0$. 
Then by the definition of $\lambda_1$ and $\Lambda_1$, $\lambda_1(Q_1, v) = (v,Q)$.

By parametrizing the oriented line $(v,Q)$ as $Q + tv$, $t \in \mathbb{R}$, we show that there exist $t_1, t_2 \in \mathbb{R}$ such that $Q + t_1v \in \text{conv}(K)^{\circ}$ and 
$Q + t_2v \in \mathbb{R}^n \setminus \overline{\text{conv}(K)}$.
This will imply that there is an intersection point of the line with $K$.

By Lemma \ref{rem:convexity}, part 2), 
$$
Q + t_1v \in \text{conv}(K)^{\circ} \iff \langle Q + t_1v, n_w \rangle > 0, \quad \forall w \in \Gamma,
$$
$$
Q + t_2v \in \mathbb{R}\setminus \overline{\text{conv}(K)} \iff \exists w \in \Gamma, \quad  \langle Q + t_2v, n_w \rangle < 0. 
$$
We define the funcion  $f(w,t) := \langle Q + tv, n_w \rangle$ for $w \in \Gamma$, $t \in \mathbb{R}$.

For $w = v$, $f(v,t) = \langle n_v, Q \rangle > 0$ is constant in $t$. For $w \neq v$, Lemma \ref{rem:convexity}, part 1) gives $\langle v, n_w \rangle > 0$, allowing us to solve the inequality $f(w,t) > 0$:
$$
\langle Q, n_w \rangle + \langle v, n_w \rangle t > 0 \quad \implies \quad t > -\frac{\langle Q, n_w \rangle}{\langle v, n_w \rangle}.
$$
We now show that there exists $t_1 >  -\frac{\langle Q, n_w \rangle}{\langle v, n_w \rangle}$ for all $w\in \Gamma\setminus \{v\}$. 

Since $\langle Q, n_v \rangle>0$ by definition of $\Delta_0$, and $\langle Q, n_w \rangle$ is continuous in $w$,
there exists a small open neighborhood $U$ of $v$ in $\Gamma$ such that $\langle Q, n_w \rangle > 0$ for all $w \in U$.
It follows that 
$$ 
-\frac{\langle Q, n_w \rangle}{\langle v, n_w \rangle}<0 
$$ 
for all $w \in U\setminus \{v\}$. 
Therefore for any $t_1 \geq 0$, we have 
$ t_1 > -\frac{\langle Q, n_w \rangle}{\langle v, n_w \rangle}$,
$f(t_1,w)>0$ for all  $w \in U\setminus \{v\}$.

For $w \in \Gamma \setminus U$, since $\langle v, n_w \rangle > 0$ and is continuous on the compact set $\Gamma \setminus U$, 
there exists a positive constant $c$ such that 
$\langle v, n_w \rangle > c$ for all $w \in \Gamma \setminus U$.
Therefore,
$$
\bigg|-\frac{\langle Q, n_w \rangle}{\langle v, n_w \rangle}\bigg| < \frac{|\langle Q, n_w \rangle|}{c} \leq \frac{\max_{w\in \Gamma}|\langle Q, n_w \rangle|}{c}.
$$
Thus, for any $t_1 \geq \frac{\max_{w\in \Gamma}|\langle Q, n_w \rangle|}{c} $, we have $t_1 > -\frac{\langle Q, n_w \rangle}{\langle v, n_w \rangle}$,
and hence 
$f(t_1,w)>0$
for all  $w \in \Gamma \setminus U$.

Therefore, take any $t_1 \geq \max\{0,\frac{\max_{w\in \Gamma}|\langle Q, n_w \rangle|}{c} \}= \frac{\max_{w\in \Gamma}|\langle Q, n_w \rangle|}{c}$, we have $f(t_1, w)>0$ for any $w\in \Gamma$.

Let us show the existence of $t_2$. 
For any $w \neq v$, we have $\langle v, n_w \rangle > 0$, which implies
$$
\lim_{t \to -\infty} f(w,t) = -\infty.
$$
Therefore, such $t_2$ must exist.

Hence, the line $Q+tv$ must intersect $K \cup \{O\}$ at some point $Q_1 = Q + t^*v, t^*\in(t_2, t_1)$.
We then have
$$
\langle n_v, Q_1 \rangle = \langle n_v, Q + t^*v \rangle = \langle n_v, Q \rangle > 0,
$$
which implies $Q_1 \neq O$. By (\ref{eq:convexity1}) in Lemma \ref{rem:convexity}, part 1), this also implies $\langle v, n_{Q_1} \rangle > 0$. 
Thus, we find $(Q_1, v) \in \Lambda_1^0$ such that $\lambda_1(Q_1, v) = (v, Q)$, establishing
$\Delta_0 \subset \lambda_1(\Lambda_1^0)$.

\vspace*{1em}
2)
From (\ref{eq:lambda11}) and (\ref{eq:lambda10}), the interior of $\Lambda_1^1$ is 
\begin{equation}\label{eq:lambda11m10}
  \Lambda_1^1\setminus \Lambda_1^0 = 
\{(Q_1,v) \mid Q_1 \in K, v\in D, \langle v, n_{Q_1} \rangle > 0 \}.
\end{equation}
To prove 2) it suffices to show that
\begin{equation}\label{eq:lambda-psi-interior}
  \lambda_1(\Lambda_1^1\setminus \Lambda_1^0) = \{(v,Q) \in TS^{n-1} \mid v \in D, Q\neq O\}.
\end{equation}
Since $v\in D\subset \text{conv}(K)^{\circ}$, 
by (\ref{eq:convexity2}) it always holds that  $\langle v, n_{Q_1} \rangle > 0$ for any $Q_1 \in K$. Therefore, (\ref{eq:lambda11m10}) can be simplified to 
$$
\Lambda_1^1\setminus \Lambda_1^0 = 
\{(Q_1,v) \mid Q_1 \in K, v\in D \}.
$$

The proof of (\ref{eq:lambda-psi-interior}) is similar to
the proof of (\ref{eq:lambda-delta}) in part 1);
in part 1), $v\in \Gamma$, while here we have $v\in D$.

First, we show ``$\subset$''. For any $(Q_1,v) \in \Lambda_1^1\setminus \Lambda_1^0$, let $(v,Q) = \lambda_1(Q_1,v)$.
By definition, 
$Q = Q_1-\langle Q_1,v\rangle v$ and $v \in D$.
It is necessary to prove 
$Q\neq O$. 
If not, then
$
Q_1= \langle Q_1,v\rangle v.
$
Taking inner product with $n_{Q_1}$ gives
$
0= \langle Q_1,v\rangle \langle v,n_{Q_1} \rangle.
$
Since $\langle v,n_{Q_1} \rangle > 0$, 
we have $\langle Q_1,v\rangle = 0$. 
Hence $Q = Q_1 \neq O$, since $Q_1 \in K$,  a contradiction.

Conversely, we show ``$\supset$''. 
For any $(v,Q) \in TS^{n-1}, v\in D, Q\neq O$, we prove that the oriented line $(v,Q)$ intersects $K$ at some point $Q_1$ where $\langle v, n_{Q_1}\rangle>0$. 

We show that there exist $t_1, t_2 \in \mathbb{R}$ such that $Q + t_1v \in \text{conv}(K)^{\circ}$ and $Q + t_2v \in \mathbb{R}^n \setminus \overline{\text{conv}(K)}$.

Now, since $v\in D$, we have $\langle v, n_w \rangle > 0$ for all $w\in \Gamma$. Thus 
solving  $f(w,t) > 0$ gives 
$$
t > -\frac{\langle Q, n_w \rangle}{\langle v, n_w \rangle}.
$$
Since $\Gamma$ is compact and the function $\big|-\frac{\langle Q, n_w \rangle}{\langle v, n_w \rangle}\big|$ is continuous on $\Gamma$, there exist $c_1$, such that 
$$
\bigg|-\frac{\langle Q, n_w \rangle}{\langle v, n_w \rangle}\bigg| < c_1.
$$ 
Thus, for any $t_1$ satisfying $t_1 \geq c_1$, we have $t_1 > -\frac{\langle Q, n_w \rangle}{\langle v, n_w \rangle}$ for all $w \in \Gamma$.
The existence of $t_2$ follows from the same argument as in part 1). Consequently, the oriented line $(v, Q)$ intersects $K$ at some point $Q_1 \in K \cup \{O\}$. To complete the proof, we must show that $Q_1 \neq O$ and $\langle v, n_{Q_1}\rangle>0$.
If $Q_1 = O$, then $\langle Q_1, v \rangle = 0$, which implies $Q_1 = Q$. This leads to $Q_1 = O$, contradicting $Q_1 \in K$.
Finally, $\langle v, n_{Q_1} \rangle > 0$ follows from the convexity of $K$ and the fact that $v$ lies inside $K$ (Lemma \ref{rem:convexity}, part 2).

\vspace*{1em}
3) We have shown in part 2) that 
$$
\psi_+^{\circ} = \{(v,Q) \in TS^{n-1} \mid v \in D, Q\neq O\},
$$
from which we can see that the boundary of $\psi_+$ 
consists of points $(v,Q)$ 
where $v$ lies on $\partial D = \Gamma$,
or $Q = O$, or both.
Thus we obtain 
\begin{equation}\label{eq:psi-plus-boundary}
  \partial \psi_+ = \{(v,Q) \in TS^{n-1} \mid v \in \Gamma\} \cup \{(v,Q) \in TS^{n-1} \mid v \in \overline{D}, Q = O\}.
\end{equation}
Let us write the first term in (\ref{eq:psi-plus-boundary}) as two parts,
$$
\begin{aligned}
  \Delta_0 &= \{(v,Q) \in TS^{n-1} \mid v \in \Gamma, \langle n_v, Q \rangle > 0\},\\
  \Delta_1 &= \{(v,Q) \in TS^{n-1} \mid v \in \Gamma, \langle n_v, Q \rangle \leq 0\}, 
\end{aligned}
$$
and denote the second term in (\ref{eq:psi-plus-boundary}) by $\Delta_{\overline{D}}$ .
Then 
$$
\partial \psi_+ = \Delta_0 \cup \Delta_1 \cup \Delta_{\overline{D}}.
$$

Lemma \ref{lem:boundary} is proved.

As a final remark of the proof, let us summarize what we have shown about the boundary of $\psi_+$.
$\Delta_0$ is the boundary between $\psi_+$ and $\psi$, and is contained in $\psi_+$;
$\Delta_1$ does not intersect $\psi_+$;
$\Delta_{\overline{D}}$
consists of all lines through the origin with directions in $\overline{D}$, 
and does not intersect $\psi_+$. The intersection properties of these sets are as follows:
$$
\Delta_0 \cap \Delta_1 = \emptyset, \quad \Delta_0 \cap \Delta_{\overline{D}} = \emptyset, \quad \Delta_1 \cap \Delta_{\overline{D}} = \{(v,Q) \in TS^{n-1} \mid v \in \Gamma , Q = O\}.
$$

\begin{remark}[The boundary structure of $\Psi$]\label{rem:boundary}
Although the boundary structure of $\Psi$ in $TS^{n-1}$ is not essential for proving integrability, 
it provides a clearer characterization of the topology of $\Psi$. 
Thus we state the following lemma without proof.
\begin{lemma}
  Let $\partial \Psi$ denote the boundary of $\Psi$ in $TS^{n-1}$.
  Then
  $$
  \partial \Psi = \Delta_1 \cup \Delta_2 \cup \Delta_3 \cup \Delta_4,
  $$
  where $\Delta_i$, $i=1,2,3,4$ are subsets of $TS^{n-1}$ defined as
  $$
  \begin{aligned}
    \Delta_1 &:= \{(v,Q) \in TS^{n-1} \mid v \in \Gamma, \langle n_v,Q \rangle \leq 0\} \cong TS^{n-2} \times \mathbb{R}_{\leq 0}, \\
    \Delta_2 &:= \{(v,Q) \in TS^{n-1} \mid v \in -\Gamma, \langle n_{-v},Q \rangle \leq 0\} \cong TS^{n-2} \times \mathbb{R}_{\leq 0}, \\
    \Delta_3 &:= \{(v,Q) \in TS^{n-1} \mid Q = O\} \cong S^{n-1}, \\
    \Delta_4 &:= \{(v,Q) \in TS^{n-1} \mid (v,Q) \text{ is tangent to } K \}. \\
  \end{aligned}
    $$
    Here, 
    $\Delta_1$ and $\Delta_2$ consist of oriented lines that are parallel to some half-line in $K$ and not intersecting $K$ transversally, with directions in $\Gamma$ and $-\Gamma$ respectively;
    $\Delta_4$ represents all oriented lines through the origin, which can be identified with the zero section of $TS^{n-1}$;
    $\Delta_3$ represents all oriented tangents to $K$.

    The subsets $\Delta_i$, $i=1,2,3$ are closed, and are submanifolds of $TS^{n-1}$. 
    Their dimensions are
  $$\dim \Delta_1 = \dim \Delta_2  = 2n-3, \quad \dim \Delta_3 = n-1.$$
  The set $\Delta_4$ is neither open nor closed. 

   The intersections among these subsets are
    $$
    \begin{aligned}
    \Delta_1 \cap \Delta_2 &= \emptyset ,\\
    \Delta_1 \cap \Delta_3 &= \{(v,Q) \in TS^{n-1} \mid v \in \Gamma, Q=O\} \cong S^{n-2}, \\
    \Delta_2 \cap \Delta_3 &= \{(v,Q) \in TS^{n-1} \mid -v \in \Gamma, Q=O\} \cong S^{n-2}, \\
    \Delta_1 \cap \Delta_4 &= \Delta_1 \cap \Delta_3, \\
    \Delta_2 \cap \Delta_4 &= \Delta_2 \cap \Delta_3, \\
    \Delta_3 \cap \Delta_4 &= (\Delta_1 \cap \Delta_3) \cup (\Delta_2 \cap \Delta_3).
    \end{aligned}
    $$
\end{lemma}
\end{remark}

\subsection{The billiard map}
In this section we study the billiard map $\mu$.
\begin{lemma}\label{lem:mu-is-diffeo}
  If the cone $K$ is $C^k$-smooth, then the billiard map $\mu: \psi_-\cup \psi \to \psi\cup \psi_+ $ is a $C^{k-1}$-smooth diffeomorphism.
\end{lemma}
\textbf{Proof:} 
  By the billiard reflection law, 
  an oriented line in $\Lambda_2$ reflects to an oriented line in $\Lambda_1$ via the map $\bar{\mu}: \Lambda_2 \to \Lambda_1$ given by
  \begin{equation}\label{eq:mu-m2-m1}
  \bar{\mu}(Q_2,v) = (Q_2, v-2\langle n_{Q_2},v \rangle n_{Q_2}).
  \end{equation}
  The map $\mu$ factors through $\bar{\mu}$ via the following commutative diagram
  $$
  \begin{tikzcd}
    \Lambda_2 \ar{r}[above]{\bar{\mu}} \ar{d}[left]{\lambda_2} & \Lambda_1 \ar{d}[right]{\lambda_1} \\
    \psi_-\cup \psi  \ar{r}[below]{\mu} & \psi_+\cup \psi 
  \end{tikzcd},
  $$
i.e.,
$$
\mu = \lambda_1 \circ \bar{\mu} \circ \lambda_2^{-1}.
$$
By Lemma \ref{lem:valid-coords}, $\lambda_1$ and $\lambda_2$ are $C^k$ diffeomorphisms. Thus it suffices to prove $\bar{\mu}$ is a $C^{k-1}$ diffeomorphism.
  
$\bar{\mu}$ is $C^{k-1}$-smooth since the normal vector field $n_{Q_2}$ is $C^{k-1}$-smooth (as $K$ is $C^k$) 
and $\bar{\mu}$ only involves smooth operations on vectors.
  
$\bar{\mu}$ is bijective with $C^{k-1}$-smooth inverse
$\nu: \Lambda_1 \to \Lambda_2$
$$
\nu(Q_1,v) = (Q_1,v-2\langle n_{Q_1},v \rangle n_{Q_1}),
$$
hence a $C^{k-1}$ diffeomorphism.

  Lemma \ref{lem:mu-is-diffeo} is proved.
\vspace*{1em}

In the next lemma, we study the billiard map in the neighborhood of $\Delta_0$.  
\begin{lemma}\label{lem:limit-of-mu}
  1)
Let $x_k$ be a sequence in $\psi$ converging to $x \in \Delta_0$ as $k\to \infty$. 
Then the limit 
$$\lim_{k\to\infty}\mu(x_k)$$ 
exists and lies in 
$\Delta_1^{\circ}=\{(v,Q) \in TS^{n-1} \mid v\in \Gamma, 
\langle n_v, Q\rangle < 0 \}\subset \Delta_1$.

2) 
The map $\sigma: \Delta_0 \to \Delta_1^\circ$ defined by the formula 
$$
\sigma(x) := \lim_{z \to x, \, z \in \psi} \mu(z),
$$ 
is well defined, and is a diffeomorphism. 

The map $\sigma$ is given by 
\begin{equation}\label{eq:sigma}
  \sigma(v,Q)=(v, Q - 2\langle Q,n_{v}\rangle  n_{v}),
\end{equation}
where $x=(v,Q)\in \Delta_0$.
\end{lemma}
\begin{figure}[htbp]
  \centering
  \includegraphics[width=0.4\linewidth]{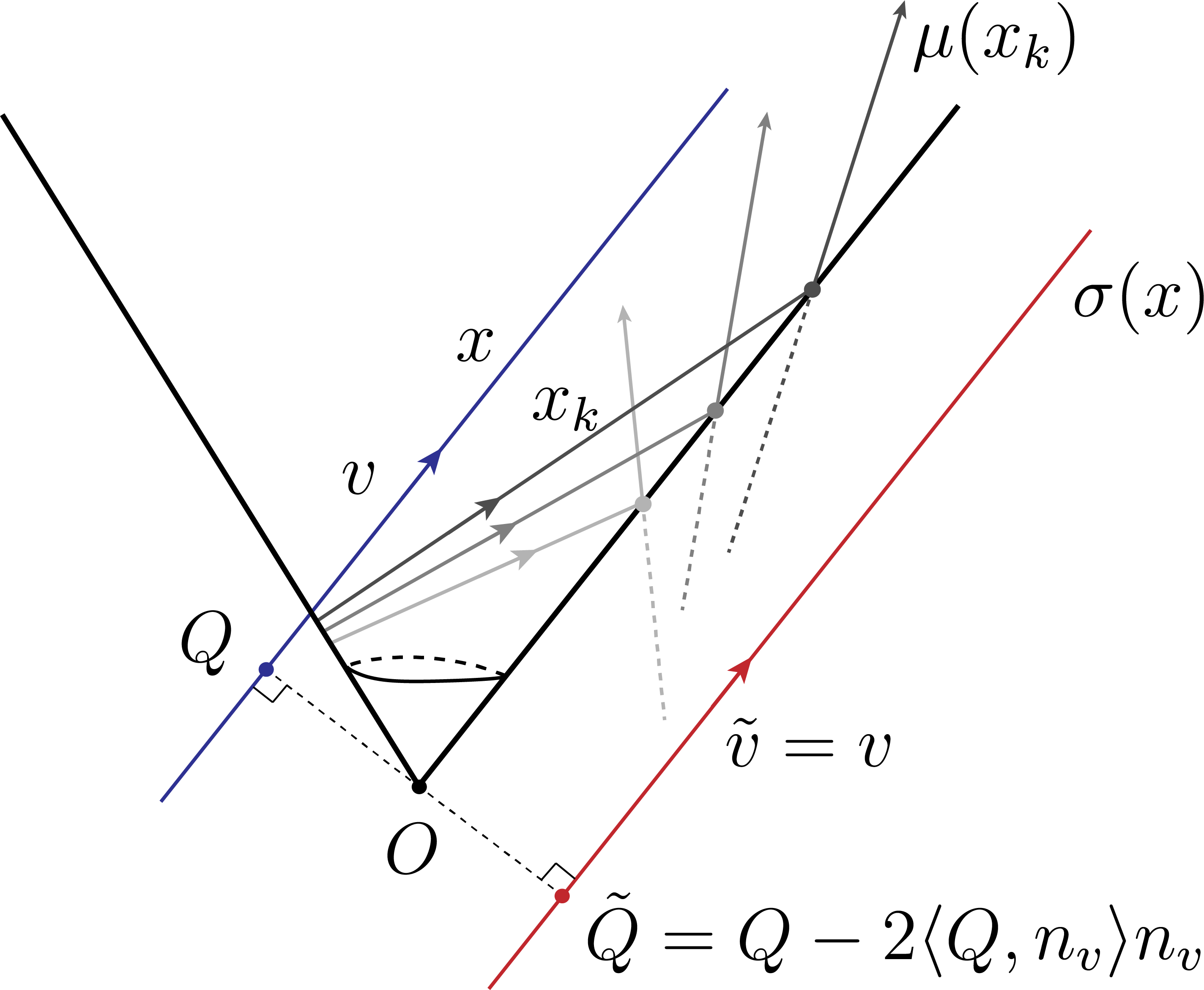}
  \caption{: If $x_k\to x \in \Delta_0$, then $\mu(x_k) \to \sigma(x)\in \Delta_1^{\circ}$.}
  \label{fig:limit-sigma}
\end{figure}
\begin{remark}
  The geometric meaning of (\ref{eq:sigma}) is the following. The line $(v,Q)$ and its image $\sigma(v,Q)$, where $v\in \Gamma$, are symmetric with respect to the tangent hyperplane to $K$ at $v$.
\end{remark}
\textbf{Proof.}
1)
Let $x_k= (v(k),Q(k)) \in \psi$ and $ (v(k),Q(k))\to x=(v,Q)\in \Delta_0$.
For each $k$, let $Q_1(k)$, $Q_2(k)$ be the intersection points of $x_k$ with $K$. 
Then $Q_1(k) \to Q_1 := K\cap x$. 
Write $Q_2(k) = t(k)w(k)$ where $t(k)>0$ and $w(k) \in \Gamma$.

We first show that $t(k) \to \infty$. If not, there exists a subsequence $t(k_j) \to t < \infty$. 
By compactness of $\Gamma$, there exists a subsequence of $k_j$ (still denoted by $k_j$) such that 
$w(k_j) \to w \in \Gamma$. Thus $Q_2(k_j) \to t w$.
Since 
$$
Q_2(k_j)-Q_1(k_j) = \|Q_2(k_j)-Q_1(k_j)\| v(k_j),
$$ 
we have 
$$
tw-Q_1 = \|tw-Q_1\|v.
$$
Note that $tw, Q_1 \in K$ and $tw, Q_1 \in x$.  
Therefore, $tw = Q_1$,  
otherwise $x$ has two intersection points with $K$,  
but this is impossible since $x \in \Delta_0$.
This implies 
$$
v=\lim_{j\to \infty} v(k_j)=\lim_{j\to \infty} \frac{Q_2(k_j)-Q_1(k_j)}{ \|Q_2(k_j)-Q_1(k_j)\|} \in T_{Q_1}K,
$$ 
hence $\langle v, n_{Q_1}\rangle=0$. However,
$\psi\cup\psi_{+}=\lambda_1(\Lambda_1)$ and by the definition of 
$\Lambda_1$, we have $\langle v, n_{Q_1}\rangle>0$, contradiction.
Hence, $t(k)\to \infty$.

Next, we prove $w(k) \to v$. From
$$
v(k)=\frac{t(k)w(k)-Q_1(k)}{\|t(k)w(k)-Q_1(k)\|},
$$
we obtain
$$
w(k)=v(k)\frac{\|t(k)w(k)-Q_1(k)\|}{t(k)} +\frac{Q_1(k)}{t(k)}.
$$
Since $\|w(k)\|=1$, the triangle inequality gives
$$
\frac{\big|t(k) - \|Q_1(k)\|\big|}{t(k)} \leq \frac{\|t(k)w(k)-Q_1(k)\|}{t(k)} \leq \frac{t(k) + \|Q_1(k)\|}{t(k)}.
$$
As $k \to \infty$, both sides converge to $1$, thus 
$\frac{\|t(k)w(k)-Q_1(k)\|}{t(k)} \to 1$. 
Combined with $\frac{Q_1(k)}{t(k)} \to 0$, we obtain $w(k) \to v$.

Let $\mu(x_k)=(\tilde{v}(k),\tilde{Q}(k))$.
We have
\begin{equation}\label{eq:tildev}
  \tilde{v}(k)= v(k)- 2\langle v(k),n_{Q_2(k)}\rangle n_{Q_2(k)}.
\end{equation}
Since $n_{Q_2(k)}= n_{w(k)}$ and $w(k)\to v$, we have 
\begin{equation}\label{eq:sigma1}
\tilde{v}(k) \to v-2\langle v,n_v\rangle n_v = v.
\end{equation}
Thus the direction of the limit of lines $x_k$ is the same as 
the limit of the direction of $\mu(x_k)$
(See Fig. \ref{fig:limit-sigma}).

The oriented line $x_k$ reflects to the oriented line $\mu(x_k)$
at the point $Q_2(k)$. On $x_k$ we have
\begin{equation}\label{eq:qk}
  Q(k) = Q_2(k)-\langle Q_2(k),{v}(k) \rangle{v}(k),
\end{equation}
(we recall that $Q(k)$ is the point of $x_k$ that realizes the distance between $x_k$ and $O$).
Similarily on $\mu(x_k)$ we have
\begin{equation}
  \label{eq:qk2}
\tilde{Q}(k) = Q_2(k)-\langle Q_2(k),\tilde{v}(k) \rangle\tilde{v}(k).
\end{equation}
Thus from (\ref{eq:qk}), (\ref{eq:qk2}) we have 
$$
\tilde{Q}(k)-  Q(k)  = \langle Q_2(k),{v}(k) \rangle{v}(k) -\langle Q_2(k),\tilde{v}(k) \rangle\tilde{v}(k).
$$
Substituting (\ref{eq:tildev}) and using $\langle Q_2(k),n_{Q_2(k)} \rangle=0$, we have
\begin{equation}
  \begin{split}
   \tilde{Q}(k)- Q(k) 
  & =   \langle Q_2(k),{v}(k) \rangle v(k) 
    -\langle Q_2(k), v(k)- 2\langle v(k),n_{Q_2(k)}\rangle n_{Q_2(k)}\rangle
   \left(v(k)- 2\langle v(k),n_{Q_2(k)}\rangle n_{Q_2(k)}\right)\\
  & =   \langle Q_2(k),{v}(k) \rangle v(k) 
    -\langle Q_2(k), v(k)\rangle
   \left(v(k)- 2\langle v(k),n_{Q_2(k)}\rangle n_{Q_2(k)}\right)\\
    \label{eq:difference-qk}
  &  =  2 \langle Q_2(k),{v}(k) \rangle  \langle v(k),n_{Q_2(k)}\rangle  n_{Q_2(k)} .
  \end{split}
\end{equation}
Taking inner product of (\ref{eq:qk}) with $n_{Q_2(k)}$ gives 
$$
\langle Q(k),n_{Q_2(k)} \rangle = - \langle Q_2(k),{v}(k) \rangle  \langle v(k),n_{Q_2(k)}\rangle.
$$
Substituting $\langle Q_2(k),{v}(k) \rangle  \langle v(k),n_{Q_2(k)}\rangle= -\langle Q(k),n_{Q_2(k)} \rangle$ into (\ref{eq:difference-qk}),
$$
\tilde{Q}(k)-  Q(k) = -2\langle Q(k),n_{Q_2(k)} \rangle  n_{Q_2(k)} ,
$$
therefore
$$
\tilde{Q}(k)=Q(k) -2\langle Q(k),n_{w(k)} \rangle  n_{w(k)} .
$$
Since $w(k)\to v$, and $Q(k)\to Q$, we have 
\begin{equation}\label{eq:sigma2}
\tilde{Q}(k)\to Q - 2\langle Q,n_{v} \rangle  n_{v} =: \tilde{Q}.
\end{equation}
From (\ref{eq:sigma1}) and (\ref{eq:sigma2}) it follows that 
$\mu(x_k)$ converges to a limit line $(v, \tilde{Q})$. 
To finish the proof of part 1) it is necessary to prove that 
$\langle n_v , \tilde{Q} \rangle < 0$.
We have
\begin{equation}\label{eq:sigma3}
  \langle n_v , \tilde{Q} \rangle = \langle n_v , Q - 2\langle Q,n_{v} \rangle  n_{v} \rangle = - \langle Q,n_{v} \rangle  <0,
\end{equation}
since $\langle Q,n_{v} \rangle>0$ by the definition of $\Delta_0$.
Hence $(v,\tilde{Q})= \lim_{k \to \infty} \mu(x_k) \in \Delta_1^{\circ}$.

\vspace*{1em}

2) We have shown by (\ref{eq:sigma1}) and (\ref{eq:sigma2}) that $\sigma$ is given by 
$$
\sigma(v,Q) =(v, Q - 2\langle Q,n_{v}\rangle  n_{v}).
$$
To show $\sigma$ is a diffeomorphism, we construct its inverse. For any $(v,Q) \in \Delta_1^{\circ}$, define
$$
\tau(v,Q) = (v, Q - 2\langle Q,n_{v}\rangle n_{v}).
$$
When $(v,Q) \in \Delta_1^{\circ}$, we have $\langle n_v,Q \rangle < 0$, thus
$$
 \langle n_v, Q - 2\langle Q,n_{v}\rangle n_{v} \rangle = -\langle n_v,Q \rangle > 0,
$$
which shows $\tau(\Delta_1^{\circ}) \subset \Delta_0$. 
A direct calculation shows that $\tau \circ \sigma = id_{\Delta_0}$ and $\sigma \circ \tau = id_{\Delta_1^{\circ}}$.
Therefore $\tau$ is the inverse of $\sigma$.
Note that $\tau$ has the same formula as $\sigma$, 
so $\sigma$ can be extended to an involution on $\Delta_0\cup \Delta_1^{\circ}\subset \partial \psi_+$.

Lemma \ref{lem:limit-of-mu} is proved.

\subsection{Construction of the first integrals}
For each $x \in \Psi$, Theorem \ref{teo:finite-intro} ensures the existence of a minimal non-negative integer $m(x)$ such that $\mu^{m(x)}(x) \in \psi_+$. 
Given a function $f: \psi_+ \to \mathbb{R}$, we define its \textit{lift to $\Psi$} as
$$
\tilde{f}(x) := f(\mu^{m(x)}(x)).
$$
Let 
$$
\Delta := \bigcup_{n=0}^{\infty} \mu^{-n}(\Delta_0).
$$ 
By Lemmas \ref{lem:boundary} and \ref{lem:mu-is-diffeo}, $\Delta$ is a codimensional-1 submanifold of $\Psi$, and thus $\Psi\setminus \Delta$ forms an open dense subset of $\Psi$. 

\begin{figure}[htbp]
  \begin{center}
  \includegraphics[scale=0.3]{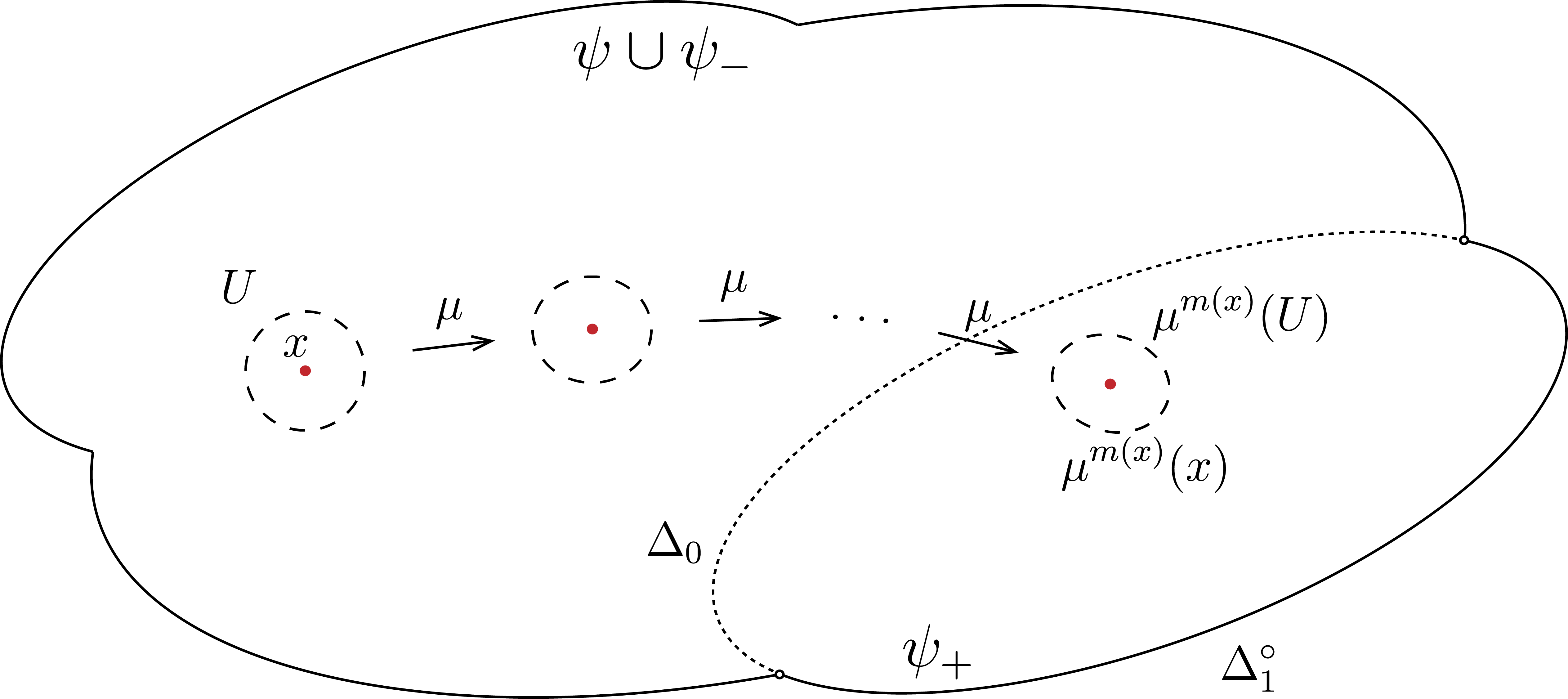}
  \end{center}
  \caption{: $x\in \Psi \setminus \Delta$, \ $\mu^{m(x)}(U)\subset \psi_+^{\circ}$.}
  \label{fig:iterations-case1}
\end{figure}

According to Lemma \ref{lem:limit-of-mu}, the lifted function $\tilde{f}$ may fail to be continuous precisely on $\Delta$. The following lemma provides a sufficient condition for the continuity of $\tilde{f}$ on $\Psi$.
Recall that by Lemma \ref{lem:boundary}
$$
\Delta_0 \subset \psi_+, \quad \Delta_1^{\circ} \subset \partial \psi_{+},
$$
and by Lemma \ref{lem:limit-of-mu}, $\sigma: \Delta_0 \to \Delta_1^{\circ}$ is a diffeomorphism.
\begin{lemma}\label{lem:lift-continuity}
  Let $f: \psi_+ \to \mathbb{R}$ be a continuous function. The function $f$ extends to $\psi_+ \cup \Delta_1^{\circ}$ via $\sigma$ by
  \begin{equation}\label{eq:extension}
  f(x) := f(\sigma^{-1}(x)), \quad  x \in \Delta_1^{\circ}.
  \end{equation}
  If $f$ is continuous on $\psi_+ \cup \Delta_1^{\circ}$, then $\tilde{f}$ is continuous on $\Psi$.
  \end{lemma}
\textbf{Proof.}
To prove $\tilde{f}$ is continuous on $\Psi$, we consider two cases:

Case 1: For $x\in \Psi \setminus \Delta$, the continuity of $\tilde{f}$ at $x$ follows directly from the continuity of $f$ and $\mu$. Indeed, since $\mu^{m(x)}(x)\in \psi_+^{\circ}$,
by continuity of $\mu$ there exists a neighborhood $U$ of $x$ such that
$\mu^{m(x)}(U)\subset \psi_+^{\circ}$
(see Fig. \ref{fig:iterations-case1}).
On $U$, $\tilde{f}$ coincides with $f\circ \mu^{m(x)}$ and is therefore continuous at $x$.

\begin{figure}[htbp]
  \begin{center}
  \includegraphics[scale=0.3]{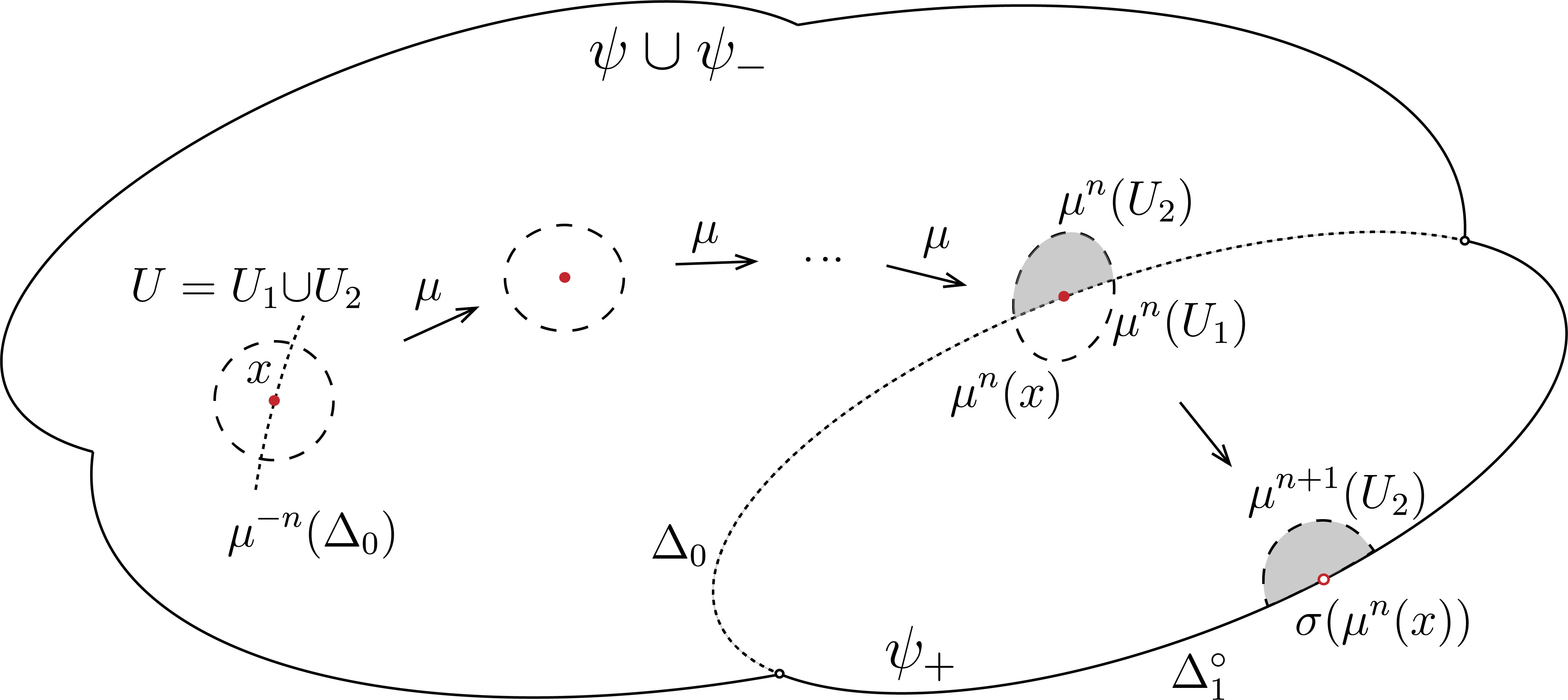}
  \end{center}
  \caption{$x\in \Delta$,\ $\mu^n(x)\in \Delta_0$; \ 
  $\mu^n(U_2)\subset \psi $, \ $\mu^n(U_1)\subset \psi_+$, \ $\mu^{n+1}(U_2)\subset \psi_+$.}
  \label{fig:iterations-case2}
\end{figure}

Case 2: For $x\in \Delta$, 
there exists a unique $n= m(x)$ such that $x \in \mu^{-n}(\Delta_0)$.
Let $U$ be a neighborhood of $x$ decomposed as
$$
U = U_1\cup U_2,
$$
where $\mu^n(U) \subset \psi_+\cup\psi$ is a neighborhood of 
$\mu^n(x)\in \Delta_0$,  $\mu^n(U_1)\subset \psi_+$, $\mu^n(U_2)\subset \psi$, 
and $\sigma(\mu^n(x))\in \Delta_1^{\circ}$ 
is a limit point of $\mu^{n+1}(U_2)\subset \psi_+$ 
(see Fig. \ref{fig:iterations-case2}).

Let $\{x_k\}$ be any sequence converging to $x$ in $U$.
We need to show that $\tilde{f}(x_k) \to \tilde{f}(x)$. 
Write $\{x_k\}$ as the union of two subsequences:
$$
\{x_k\} = \{x_{k_1}\} \cup \{x_{k_2}\},
$$
where $x_{k_1} \in U_1$ and $x_{k_2} \in U_2$.

There are three cases to consider.

\begin{figure}[htbp]
	\centering
	\begin{minipage}[t]{0.48\linewidth}
		\centering
		\includegraphics[width=1\textwidth]{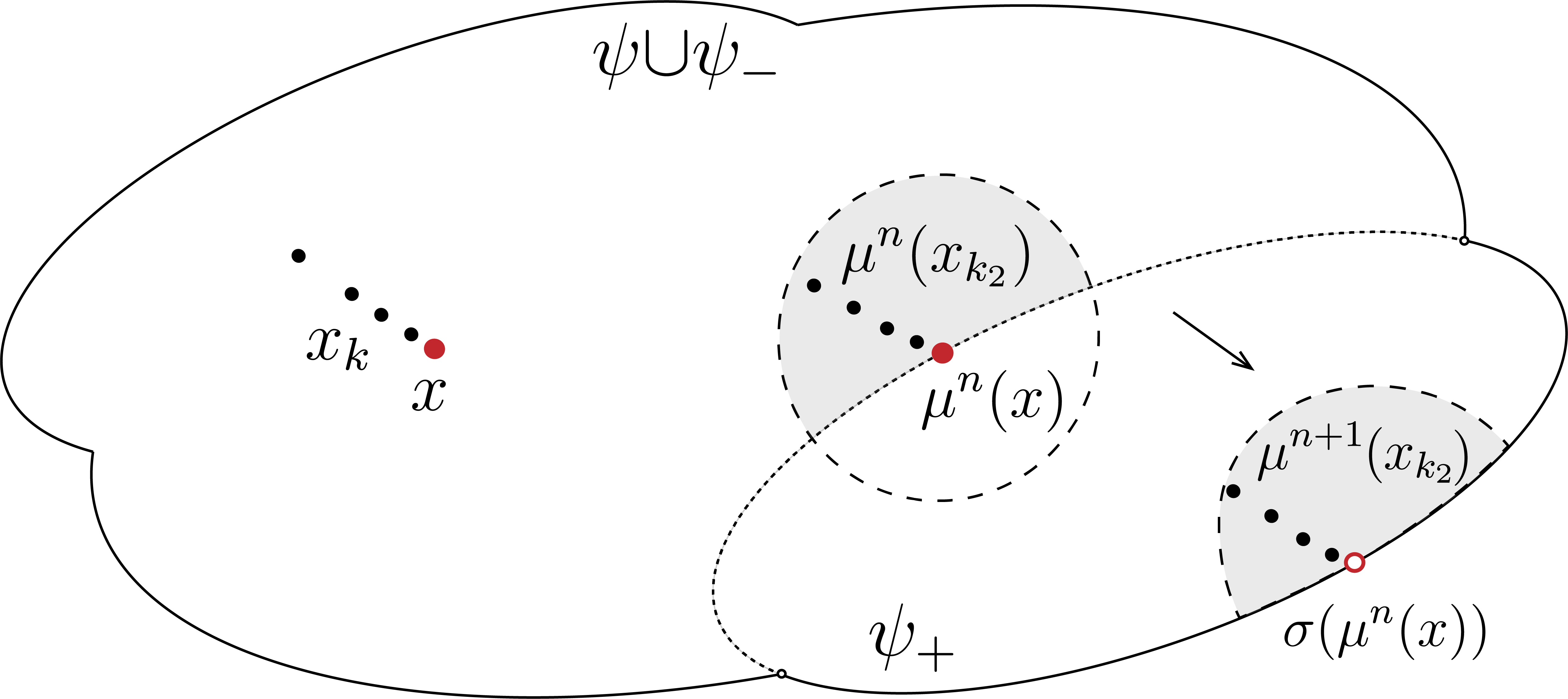}
    \caption{(a): $\mu^n(x)\in \psi_+$, $\mu^{n+1}(x_{k_2})\in \psi_+$.}
    \label{fig:iterations-a}
	\end{minipage}
  \begin{minipage}[t]{0.48\linewidth}
		\centering
		\includegraphics[width=1\textwidth]{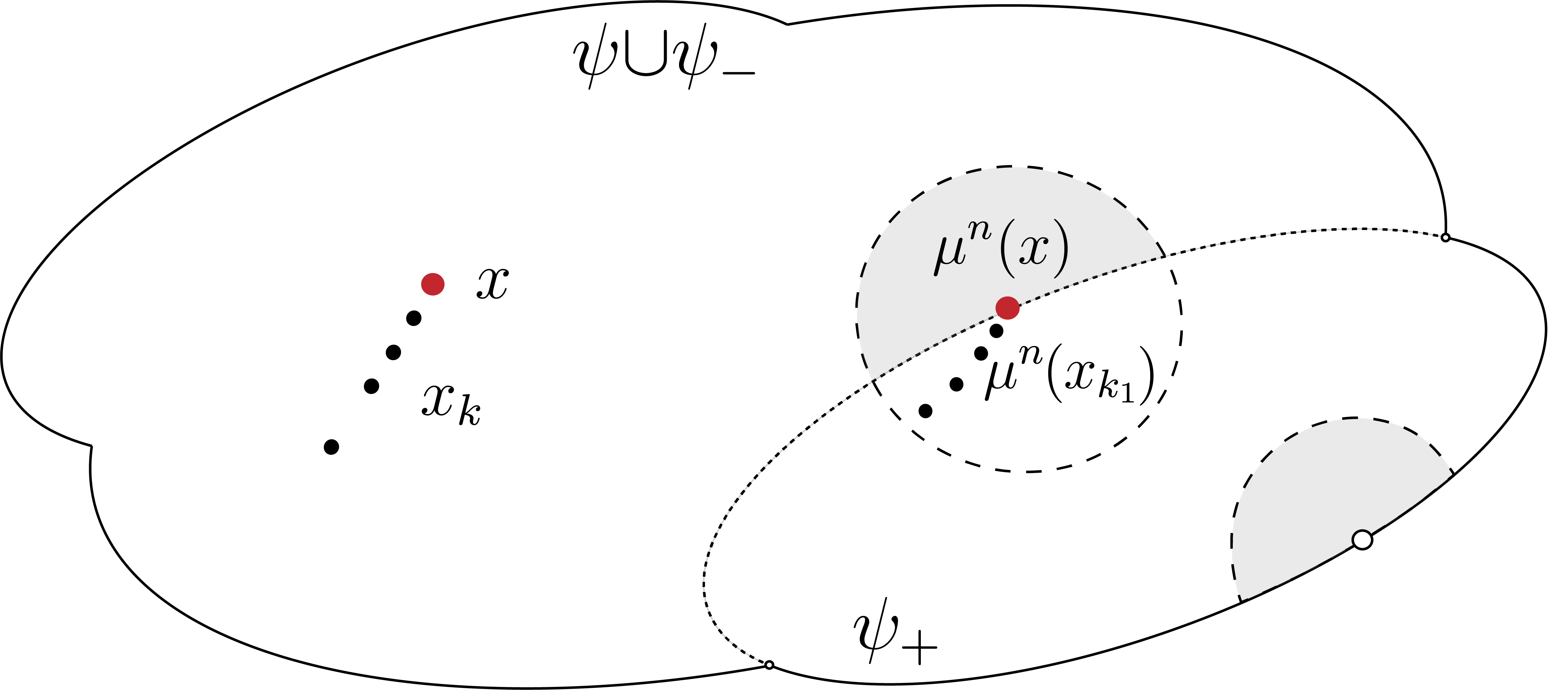}
    \caption{(b): $\mu^n(x)\in \psi_+$, $\mu^{n}(x_{k_1})\in \psi_+$.}
    \label{fig:iterations-b}
	\end{minipage}
\end{figure}

(a) 
If $\{x_{k_1}\}$ is finite, 
then there exists $k_0$ such that $x_k \in U_2$ for all $k \geq k_0$.
We have $\mu^{n+1}(x_{k_2})\in \psi_+$ (see Fig. \ref{fig:iterations-a}).
By definition, the value of $\tilde{f}$ at each $x_{k_2}$ is 
\begin{equation}\label{eq:f-tilde-k1}
\tilde{f}(x_{k_2}) = f(\mu^{n+1}(x_{k_2})),
\end{equation}
while the value at $x$ is
\begin{equation}\label{eq:f-tilde-x}
\tilde{f}(x) = f(\mu^{n}(x)).
\end{equation}
Since $\mu^{n}(x_{k_2})\in \psi$ and $\mu^{n}(x_{k_2}) \to \mu^{n}(x) \in \Delta_{0}$ (see Fig. \ref{fig:iterations-a}), by the definition of $\sigma$ 
(see Lemma \ref{lem:limit-of-mu}), 
we have
\begin{equation}\label{eq:mu-limit}
\mu^{n+1}(x_{k_2}) \to \sigma(\mu^{n}(x)).
\end{equation}
Then, by  (\ref{eq:f-tilde-k1}), (\ref{eq:mu-limit}), and the continuity of $f$ at $\sigma(\mu^{n}(x))$, we have
\begin{equation}\label{eq:f-tilde-limit}
\lim_{k_2\to\infty}\tilde{f}(x_{k_2}) 
= \lim_{k_2\to\infty}f(\mu^{n+1}(x_{k_2})) = f(\sigma(\mu^n(x))).
\end{equation}
By the extension rule (\ref{eq:extension}) of $f$, we have 
\begin{equation}\label{eq:extension-rule}
f(\sigma(\mu^n(x))) = f(\mu^n(x)).
\end{equation}
Thus, combining (\ref{eq:extension-rule}), 
(\ref{eq:f-tilde-limit}), (\ref{eq:f-tilde-x}),  
we have 
$\lim_{k_2\to\infty}\tilde{f}(x_{k_2}) = \tilde{f}(x)$
, i.e., 
$\tilde{f}$ is continuous at $x$.

(b) 
If $\{x_{k_2}\}$ is finite, then there exists $k_0$ such that $x_k \in U_1$ for all $k \geq k_0$. See Fig. \ref{fig:iterations-b}.
We have 
\begin{equation}\label{eq:direct-limit}
\lim_{k\to\infty}\tilde{f}(x_k) = \lim_{k_1\to\infty}f(\mu^{n}(x_{k_1})) = f(\mu^{n}(x)) = \tilde{f}(x),
\end{equation}
as $x_{k_1}$ and $x$ are mapped by the same $n$ iterations to $\psi_+$.

(c) 
If both subsequences are infinite, then from (\ref{eq:f-tilde-limit}) we have 
\begin{equation}\label{eq:limit-seq1}
\tilde{f}(x_{k_2}) \to f(\sigma(\mu^n(x))),
\end{equation}
and from (\ref{eq:direct-limit}) we have 
\begin{equation}\label{eq:limit-seq2}
\tilde{f}(x_{k_1}) \to f(\mu^n(x))= \tilde{f}(x).
\end{equation}
By the extension rule  (\ref{eq:extension})  of $f$ on $\psi_+ \cup \Delta_1^{\circ}$,
these two limits are equal:
\begin{equation}\label{eq:equal-limits}
f(\mu^n(x)) = f(\sigma(\mu^n(x)))= \tilde{f}(x).
\end{equation}
Therefore, $\tilde{f}$ is continuous at $x$.

Lemma \ref{lem:lift-continuity} is proved.

\begin{remark}\label{rmk:smooth}
  If $f$ is smooth on $\psi_+^{\circ}$, then $\tilde{f}$ is smooth on $\Psi\setminus \Delta$. This follows from the proof above: for $x \in \Psi \setminus \Delta$, $\tilde{f}$ locally coincides with the composition of smooth functions $f\circ \mu^{m(x)}$.
\end{remark}

On the boundary $\Delta_0$, 
the extension condition in 
Lemma \ref{lem:lift-continuity}, combined with (\ref{eq:sigma}),
requires that
\begin{equation}\label{eq:glue}
f(v,Q) = f(v,Q - 2\langle Q,n_v\rangle n_v)
\end{equation}
for all $(v,Q) \in \Delta_0$, 
i.e., if $f$ satisfies (\ref{eq:glue}), then $\tilde{f}$ is continuous. 
Let us construct examples of the continuous first integrals.

\begin{example}
  Let us check that $I=\sum_{i<j}m^2_{i,j}$ 
  (see Theorem \ref{teo:cone-intro-n})
  satisfies the condition (\ref{eq:glue}).
Let $(v,Q)\in TS^{n-1}$. Then 
$$
I(v,Q)= \langle Q, Q \rangle.
$$
Let $(v,Q) \in \Delta_0$, then 
$$
I(\sigma(v,Q)) = \langle Q - 2\langle Q, n_v \rangle n_v, Q - 2\langle Q, n_v \rangle n_v \rangle
=  \langle Q, Q \rangle = I(v,Q).
$$
\end{example}

\begin{example}  
For any trajectory in the cone,
let $v=(v^1,\dots,v^{n})$ denote the direction of its final oriented line. 
The function $v^i: \psi_+ \to \mathbb{R}$ can be naturally extended to the function $v^i: \psi _+ \cup \Delta_1^{\circ} \to \mathbb{R}$, and this function satisfies  (\ref{eq:glue}). 
Hence, by Lemma \ref{lem:lift-continuity}, we have continuous first intergral $\tilde{v}^i$ on $\Psi$. 
Thus, any continuous function of $v$ defines a continuous
first intergral. 
\end{example}

\begin{example} 
  Let us define functions on $\Delta_0 \cup \Delta_1^{\circ}$:
  \begin{equation}\label{eq:f_i}
    f_i(v,Q) = (Q-\langle Q,n_v\rangle n_v)^i,\quad i=1,\ldots,n,
    \quad (v,Q)\in \Delta_0\cup \Delta_1^{\circ}.
  \end{equation}
  Here, 
  $$
  (f_1,\ldots, f_n)= Q-\langle Q,n_v\rangle n_v
  $$
  is the projection of $Q$ on $T_vK$. 
  
  The functions $f_i$, $i=1,\ldots,n$, satisfy (\ref{eq:glue}).
  Indeed,
  $$
  \begin{aligned}
    f_i(v,Q-2\langle Q,n_v\rangle n_v) 
  & = 
    \left( Q-2\langle Q,n_v\rangle n_v -\langle Q-2\langle Q,n_v\rangle n_v ,n_v\rangle n_v \right)^i \\
  & = \left( Q-2\langle Q,n_v\rangle n_v - \langle Q,n_v\rangle n_v + 2\langle Q,n_v\rangle n_v \right)^i \\
  & = f_i(v,Q).
  \end{aligned}
  $$
  We will use $f_i(v,Q)$ to construct continuous first integrals.
  \end{example}

\begin{figure}[htbp]
  \begin{center}
  \includegraphics[scale=0.3]{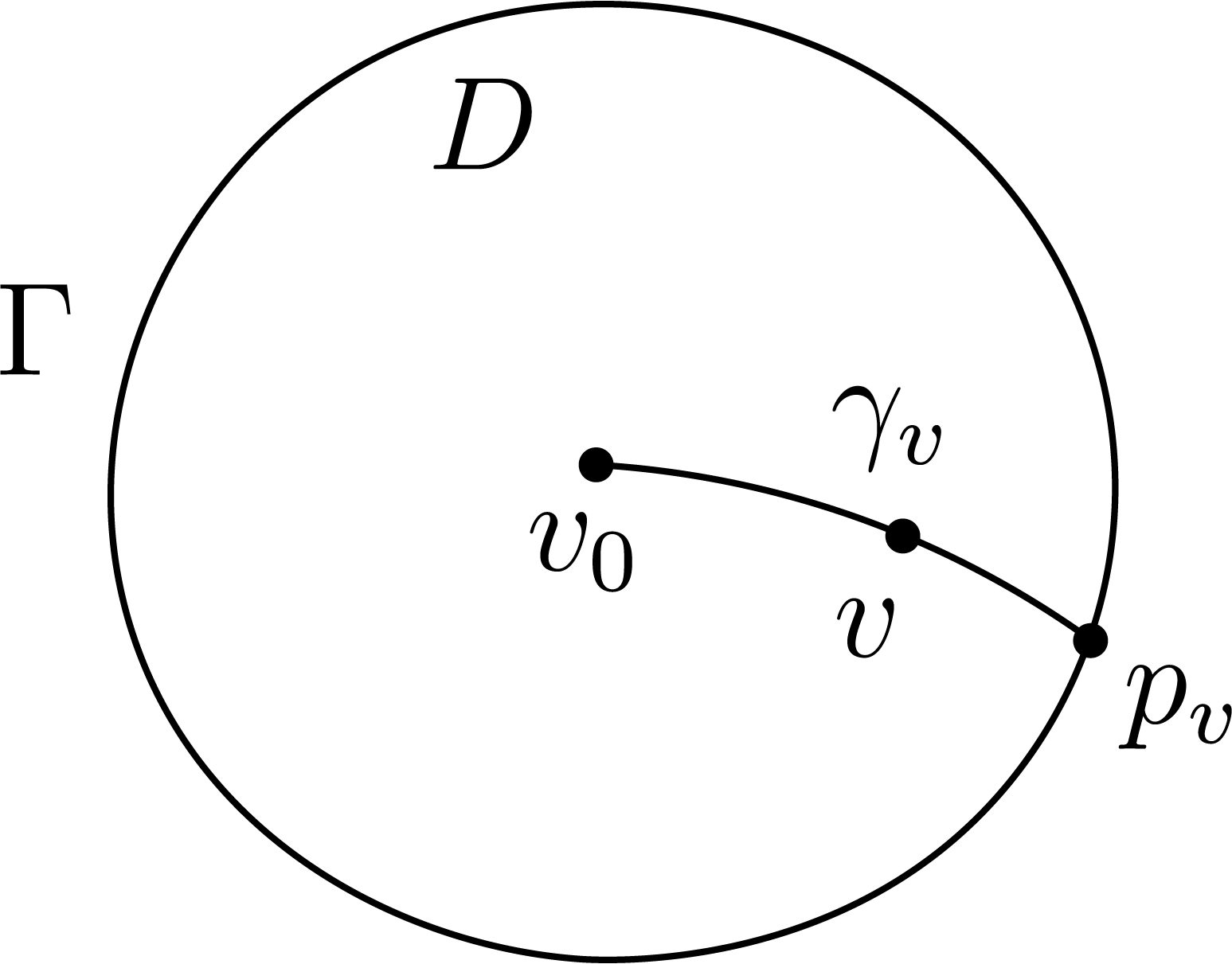}
  \end{center}
  \caption{The geodesic $\gamma_v$.}
  \label{fig:geodesic}
\end{figure}

Next, we extend the functions $f_i$ in Example 3 continuously to $\psi_+$, by constructing a continuous tangent vector field on $\overline{D}$ extending $n_v$ from $\Gamma$. By the Poincar\'e-Hopf theorem, since $\overline{D}$ has Euler characteristic 1, any such extension must have at least one vanishing point. We construct such an extension as follows:

Fix $v_0=(0,\ldots,0,1)\in D \subset S^{n-1}$. 
For any $v \in D$, let $\gamma_v$ be the unique minimal geodesic in $S^{n-1}$ that starts
 from $v_0$, passing through $v$, and intersects $\Gamma$ at  $p_v$ (see Fig. \ref{fig:geodesic}). When $v\in \Gamma$, $v=p_v$.

Define
\begin{equation}\label{eq:distance-on-d}
  g(v) = \frac{d(v_0,v)}{d(v_0,p_v)},
\end{equation}
where $d(\cdot,\cdot)$ denotes the geodesic distance on $S^{n-1}$. Since $\Gamma$ is $C^3$-smooth, the function $g: \overline{D} \to [0,1]$ is $C^3$-smooth. 
Note that $g(v_0)=0$, $g(v)=1$ for $v \in \Gamma$, and $g$ takes values in $[0,1)$ on $D$.

For any $v\in \overline{D}$, define $X(v)$ to be the vector obtained by parallel transporting the unit normal vector $n_{p_v}$ along $\gamma_v$ to $v$
multiplied by $g(v)$. 
Note that  $\|X(v)\| = g(v)$ for $v \in \overline{D}$, and $X(v)= n_v$ for $v \in \Gamma$. 

On $\psi_+ \cup \Delta_1^{\circ}$, using the vector field $X(v), v\in \overline{D}$, we define
\begin{equation}\label{eq:construction-symmetry}
  F_i(v,Q) = (Q-\langle Q,X(v)\rangle X(v))^i,\quad i=1,\ldots,n.
\end{equation}
These $F_i$ provide continuous extensions of $f_i$ to 
$\psi_+ \cup \Delta_1^{\circ}$.  
By Example 3 and Lemma \ref{lem:lift-continuity}, functions
$F_i$ define continuous integrals $\tilde{F}_i$ on $\Psi$.
In fact, as the construction of $X(v)$ is smooth in $v\in D$, 
the functions $F_i$ are smooth on $\psi_{+}^{\circ}$ and their lifts 
$\tilde{F}_i$ are smooth on $\Psi\setminus \Delta$ (see Remark \ref{rmk:smooth}).

\begin{lemma}\label{lem:continuous1}
For fixed $v\in D$, the map $(F_1,\ldots,F_{n-1})$ is an invertible linear transformation from $T_vS^{n-1}$ to $\mathbb{R}^{n-1}$. 
For fixed $v\in \Gamma$, the map $(F_1,\ldots,F_{n-1})$ is a  linear transformation from $T_vS^{n-1}$ to $\mathbb{R}^{n-1}$ of 
rank $n-2$.
\end{lemma}
\textbf{Proof.} 
For any fixed $v\in \overline{D}$, $T_vS^{n-1}$ is the $n-1$ dimensional linear subspace of $\mathbb{R}^n$ defined by 
$$
\langle v,Q \rangle =0.
$$
Since $v^n>0$ for  $v\in \overline{D}$, we have 
\begin{equation}\label{eq:qn-v}
  Q^n = -\frac{1}{v^n}\sum_{i=1}^{n-1} Q^iv^i.
\end{equation}
The map $(F_1,\ldots,F_{n-1})$ is linear in $Q$ by definition.

To determine the rank of the linear transformation $(F_1,\ldots,F_{n-1})(Q): T_vS^{n-1} \to \mathbb{R}^{n-1}$, we consider three cases: $v = v_0$, $v \in D \setminus \{v_0\}$, and $v \in \Gamma$.

For any vector $Y=(Y^1,\ldots,Y^{n})\in \mathbb{R}^n$, let $\hat{Y}=(Y^1,\ldots,Y^{n-1})$ denote the vector consisting of the first $n-1$ components of $Y$.

1. At $v=v_0$ where $X(v)=0$, we have 
$$(F_1,\ldots,F_{n-1})(Q)=\hat{Q}.$$
If $\hat{Q}=0$, then by (\ref{eq:qn-v}) we have $Q^n=0$. Then 
$Q=0$, which implies 
$$
\ker (F_1,\ldots,F_{n-1}) = \{0\}.
$$
Therefore,
$$
\text{rank} (F_1,\ldots,F_{n-1})= n-1.
$$

2. At $v\in D\setminus \{v_0\}$, $0<\|X(v)\|<1$.
Decompose $T_v S^{n-1}$ into the direct sum $V_1\oplus V_2$, where $V_1$ is the one-dimensional subspace spanned by $X(v)$ and $V_2$ is its orthogonal complement.
We now show that if $(F_1,\ldots,F_{n-1})(Q)=0$, then $Q=0$. 

For any $Q\in T_vS^{n-1}$, write $Q= Q^{(1)}+Q^{(2)}$,
where $Q^{(1)}\in V_1$ and $Q^{(2)}\in V_2$.

Since
$$\langle Q^{(1)},X(v)\rangle X(v) = \|X(v)\|^2 Q^{(1)},
\quad 
\langle Q^{(2)},X(v)\rangle=0,$$ 
we have 
\begin{equation}\label{eq:f-q1-q2}
  (F_1,\ldots,F_{n-1})(Q) = {(1-\|X(v)\|^2)\hat{Q}^{(1)}+\hat{Q}^{(2)}}.
\end{equation}
If  $(F_1,\ldots,F_{n-1})(Q)=0$ then ${(1-\|X(v)\|^2)\hat{Q}^{(1)}+\hat{Q}^{(2)}}=0$, 
and then by (\ref{eq:qn-v})  
\begin{equation}\label{eq:decom-q1-q2}
  {(1-\|X(v)\|^2)Q^{(1)}+Q^{(2)}}=0.
\end{equation}
Since $1-\|X(v)\|^2 \neq 0$, $Q^{(1)}\in V_1$, $Q^{(2)}\in V_2$, from (\ref{eq:decom-q1-q2})
we have $Q^{(1)}=0$ and $Q^{(2)}=0$. 
Therefore
$$
\ker (F_1,\ldots,F_{n-1}) = \{0\}.
$$
Therefore,
$$
\text{rank} (F_1,\ldots,F_{n-1})= n-1.
$$

  3. At $v\in \Gamma$, since $\|X(v)\|=\|n_v\|=1$, 
  by (\ref{eq:f-q1-q2}) we have 
  $$
  (F_1,\ldots,F_{n-1})(Q) = (1-\|X(v)\|^2)\hat{Q}^{(1)}+\hat{Q}^{(2)}= \hat{Q}^{(2)}.
  $$
 If  $(F_1,\ldots,F_{n-1})(Q)=0$ then  $\hat{Q}^{(2)}=0$, and then from (\ref{eq:qn-v}) we have 
 $Q^{(2)}=0$. Therefore,
 $$
 \ker (F_1,\ldots,F_{n-1}) = V_1.
 $$
  Therefore,
  $$
\text{rank} (F_1,\ldots,F_{n-1})= n-2.
$$

  Lemma \ref{lem:continuous1} is proved.

  \vspace*{1em}

By Lemma \ref{lem:continuous1}, the functions 
$F_1,\ldots,F_{n-1}$ constructed in 
(\ref{eq:construction-symmetry}) are functionally independent for fixed $v \in D$.
Define functions $I_i: \Psi \to \mathbb{R}$, $i=1,\ldots,2n-1$ as follows:
$$
I_i(v,Q) = \tilde{v}^i, \quad i = 1, \ldots, n-1,
$$
where $\tilde{v}^i$ are constructed in Example 2;
$$
I_i(v,Q) = \tilde{F}_{i-(n-1)}(v,Q), \quad i = n, \ldots, 2n-2,
$$
where $\tilde{F}_j$ are the lifts of $F_j$ from $\psi_+$ to $\Psi$;
and
$$
I_{2n-1}(v,Q) = \langle Q, Q \rangle,
$$
which coincides with the first integral $I=\sum_{i<j}m^2_{i,j}$ in Theorem \ref{teo:cone-intro-n}.

\begin{lemma}\label{lem:complete-integrals}
The functions $I_1,\ldots,I_{2n-1}$ defined above are continuous first integrals on $\Psi$. Moreover, 
the map 
$$
\mathcal{I}:=(I_1,\ldots,I_{2n-1}): \Psi \setminus \Delta \to \mathbb{R}^{2n-1}
$$ 
is a smooth immersion, 
and for any  $c_1,\ldots, c_{2n-1}\in \mathbb{R}$, the level set
$$\{x \in \Psi \mid I_i(x)=c_i, i=1,\ldots,2n-1\}$$
is either empty or consists of exactly one billiard trajectory.
\end{lemma}
\textbf{Proof.} 
The continuity and $\mu$-invariance of $I_1,\ldots,I_{2n-1}$ on $\Psi$ follow from Example 2, Example 3, 
the construction of $F_j$, and Theorem 1.

Let us prove that $\mathcal{I}$ is an immersion. 
Write $v=(v^1,\ldots,v^n)$ and $Q=(Q^1,\ldots,Q^n)$. The submanifold 
$TS^{n-1}\subset \mathbb{R}^{n}\times \mathbb{R}^n$ is defined by
$$
\langle v,v \rangle = 1, 
\quad  
\langle v,Q \rangle =0.
$$
Since $v^n>0$ on $\overline{D}$, we have
$$
v^n = \sqrt{1-\sum_{i=1}^{n-1}(v^i)^2},
\quad
Q^n = -\frac{1}{v^n}\sum_{i=1}^{n-1} Q^iv^i.
$$
Therefore, $(v^1,\ldots,v^{n-1},Q^1,\ldots,Q^{n-1})$ forms a local coordinate system 
on $\psi_+ \cup \Delta_1^{\circ}$.

By definition $v^i, F_j, 1\leq i,j\leq n-1$ are smooth functions on 
$\psi_{+}^{\circ}$.
The Jacobian matrix of $(v^1,\ldots v^{n-1}, F_1, \ldots, F_{n-1})$ with respect to  $(v^1,\ldots, v^{n-1}, Q^1, \ldots, Q^{n-1})$ is 
$$
\begin{pmatrix}
\frac{\partial v^i}{\partial v^k}
&
\frac{\partial v^i}{\partial Q^l} \\
\frac{\partial F_j}{\partial v^k}
&
\frac{\partial F_j}{\partial Q^l} 
\end{pmatrix}
=
\begin{pmatrix}
  I_{n-1} & 0 \\
  * &  \frac{\partial F_j}{\partial Q^l}
\end{pmatrix}.
$$
By Lemma \ref{lem:continuous1}, we know that  
$$
\text{rank} \left(\frac{\partial F_j}{\partial Q^l}\right)=n-1
$$
for all $(v,Q)\in \psi_{+}^{\circ}$. Thus $v^1,\ldots v^{n-1}, F_1, \ldots, F_{n-1}$ are functionally independent on $\psi_{+}^{\circ}$, and hence 
$I_1,\ldots,I_{2n-2}$ 
are functionally independent on $\Psi\setminus \Delta$.
Therefore the map $\mathcal{I}=(I_1,\ldots,I_{2n-1})$ is a smooth immersion on $\Psi\setminus \Delta$.

It remains to prove that the value of $\mathcal{I}$
uniquely determines a trajectory. 
Equivalently, the equations  
$$  
\begin{cases}
v^i = c_i  & i=1,\ldots,n-1 \\
F_i(v,Q) = c_{n+i-1}  & i=1,\ldots,n-1 \\
\langle Q,Q \rangle = c_{2n-1}
\end{cases}
$$
uniquely determine $(v, Q) \in \psi_+$.

First, $v^1=c_1,\ldots,v^{n-1}=c_{n-1}$ uniquely determine $v\in D\cup\Gamma$. 

When $v\in D$, by Lemma \ref{lem:continuous1}, $F_1,\ldots,F_{n-1}$ define an invertible linear map from 
$T_vS^{n-1}$ to $\mathbb{R}^{n-1}$ at each $v$. Thus the equations
$$
F_i(v,Q) = c_{n+i-1}, \quad i=1,\ldots,n-1,
$$
uniquely determine $Q \in T_vS^{n-1}$.

When $v\in\Gamma$, i.e., $(v,Q)\in \Delta_0$, the functions $F_1,\ldots,F_{n}$ coincide with 
$f_1,\ldots,f_{n}$ in Example 3, which are the coordinates of $Q-\langle Q,n_v\rangle n_v\in T_vK$ (see Example 3).
Since 
$\langle Q-\langle Q,n_v\rangle n_v, v \rangle = 0$, we have 
$$
\sum_{i=1}^{n} v^i F_i = 0,
$$
thus 
$$
F_n = -\frac{1}{v^n}\sum_{i=1}^{n-1}v^iF_{i}. 
$$
Therefore,
the values of $F_1,\ldots,F_{n-1}$ determine $F_n$ uniquely, and hence uniquely determine 
\begin{equation}\label{eq:qminus}
  Q-\langle Q,n_v\rangle n_v.
\end{equation}
Since $(v,Q) \in \Delta_0$, we have $\langle Q,n_v\rangle>0$ (see Lemma \ref{lem:boundary}). 
Then the value of $I_{2n-1}=\langle Q,Q\rangle$ determines $\langle Q,n_v\rangle$ via
$$\langle Q,n_v\rangle = \sqrt{\langle Q,Q\rangle - \|Q-\langle Q,n_v\rangle n_v\|^2}.$$
From (\ref{eq:qminus}), $Q$ is uniquely determined.

Lemma \ref{lem:complete-integrals} is proved.

\vspace*{1em}

This proves Theorem \ref{teo:integrability}.

\vspace*{1em}

In the next lemma, we construct functions on $\psi_+$ that lift to smooth integrals on $\Psi$.
Through a suitable auxiliary function $h$, we obtain $2n-1$ integrals $(I^s_1,\ldots,I^s_{2n-1})$ such that a unique trajectory in $\Psi\setminus\Delta$ is determined by the value of 
$(I^s_1,\ldots,I^s_{2n-1}) \in \mathbb{R}^{2n-1}$.

To define $h$, we begin with a function $g_1: \overline{D} \to [0,1]$ given by
\begin{equation}\label{eq:geo-pv}
g_1(v) = \frac{\sqrt{(v^1)^2+\cdots+(v^{n-1})^2}}{\sqrt{(p_v^1)^2+\cdots+(p_v^{n-1})^2}},
\end{equation}
where for each $v \in \overline{D}$, 
$p_v$ is the intersection point of $\Gamma$ with the extension of the unique minimal geodesic in $S^{n-1}$ from $v_0=(0,\ldots,0,1)$ to $v$ (see Fig. \ref{fig:geodesic}).
The function $g_1$ is $C^3$-smooth on $\overline{D}$ and has properties similar to $g$: 
it equals $1$ on $\Gamma$, equals $0$ at $v_0$, and takes values in $[0,1)$ on $D$.

Let $h(v) := (g_1(v)-1)^2$. Define functions $v^i_h,  Q^i_h:\psi_{+} \to \mathbb{R}$ by 
$$
v^i_h(v,Q) = v^i h(v), \quad
Q^i_h(v,Q) =Q^{i}h(v), \quad i = 1,\ldots,n-1.
$$
For $i=1,\ldots,2n-1$, let $I^s_i: \Psi \to \mathbb{R}$ be defined by
\begin{equation*}\label{eq:intergrals3}
\begin{aligned}
&I^s_i(v,Q) = \tilde{v}^i_h, \quad i = 1,\ldots,n-1,\\
&I^s_i(v,Q) = \tilde{Q}^{i-(n-1)}_h, \quad i = n,\ldots,2n-2,\\
&I^s_{2n-1}(v,Q)  = \tilde{h},
\end{aligned}
\end{equation*}
where $\tilde{v}^i_h$, $\tilde{Q}^j_h$ and $\tilde{h}$ are the lifts of $v^i_h$, $Q^j_h$ and $h$ to $\Psi$, respectively.
\begin{lemma}\label{lem:many-integrals}
  The functions $I^s_1,\ldots,I^s_{2n-1}$ defined above are $C^1$-smooth first integrals on $\Psi$, which vanish identically on $\Delta$. Moreover, for any $c_1,\ldots, c_{2n-1}\in \mathbb{R}$, the level set
  $$
  \{x \in \Psi\setminus \Delta \mid I^s_i(x)=c_i, i=1,\ldots,2n-1\}
  $$
  is either empty or consists of exactly one billiard trajectory.
\end{lemma}
\textbf{Proof.} 
  For any $x=(v,Q)\in \Delta_0\cup \Delta_1^{\circ}$, 
  $h(v)=0$ and $dh(v)=0$,
  and therefore
  $$
  v^ih(v)=0,\quad Q^ih(v)=0.
  $$
  Hence by Lemma \ref{lem:lift-continuity}, $I_i^s$ are continuous on $\Psi$.
  We have on $\Delta_0\cup \Delta_1^{\circ}$ 
$$  
d(v^i h)= h dv^i + v^i dh =0, 
\quad d(Q^i h)= h dQ^i + Q^i dh =0,
$$ 
hence $I_i^s$ are $C^1$ smooth on $\Psi$.
Therefore, 
$I^s_1,\ldots,I^s_{2n-1}$ are $C^1$-smooth first integrals on $\Psi$. They vanish on $\Delta$, as 
$h$ vanishes on $\Delta_0$ by definition.

To determine the level set, for $(v,Q)\in \psi_{+}^{\circ}$, 
we need to solve
$$
v^ih(v)= c_i, \quad 
Q^i h(v) = c_{n-1+i}, \quad 
h(v) = c_{2n-1}, 
\quad
i=1,\ldots,n-1.
$$
If $c_{2n-1}=0$ the system has no solution in $\psi_{+}^{\circ}$ 
as $h(v)$ is non-zero on $\psi_{+}^{\circ}$.
If $c_{2n-1}\neq 0$, we obtain
$$
v^i=\frac{c_i}{c_{2n-1}}, \quad
Q^i = \frac{c_{n-1+i}}{c_{2n-1}}, \quad i=1,\ldots,n-1,
$$
which corresponds to at most a unique point in $\psi_{+}^{\circ}$.
Therefore, the level set is either empty or consists of exactly one billiard trajectory.

Lemma \ref{lem:many-integrals} is proved.

\vspace{2cm}

\noindent
{Andrey E. Mironov}\\
Sobolev Institute of Mathematics, Novosibirsk, Russia\\
Email: \texttt{mironov@math.nsc.ru}

\medskip
\noindent
{Siyao Yin}\\
Sobolev Institute of Mathematics, Novosibirsk, Russia\\
Email: \texttt{yinsiyao@outlook.com}

\end{document}